\documentclass[11pt]{amsart}
\usepackage{amssymb,mathrsfs,amsmath}
 \usepackage{graphicx}
\usepackage{amscd,amsmath,amsopn,amssymb,amsthm,multicol}
\usepackage[color,matrix, all, 2cell]{xy}
\usepackage{amscd}
\usepackage{lscape}
\usepackage{slashed}
\usepackage{graphicx}
\usepackage{setspace}
\usepackage{upgreek}
\usepackage{textgreek}
\usepackage{enumerate}
\usepackage{color}
\usepackage{lscape}
\usepackage{tikz}
\usepackage{multirow}
\usepackage{cancel}
\usepackage{soul}
\usepackage{comment}
\usepackage{wasysym}
\usepackage{mathrsfs}
\usepackage{young}
\usepackage{mathtools}
\usepackage{bbm}
\usepackage{marginnote}
\usepackage{diagbox}

\numberwithin{equation}{section}

 \newcommand{\R}{\mathbb{R}}
\newcommand{\Z}{\mathbb{Z}}
\newcommand{\C}{\mathbb{C}}
\newcommand{\Hn}{\mathbb{H}}

\renewcommand{\so}{\mathfrak{so}}
\renewcommand{\sp}{\mathfrak{sp}}

\DeclareMathOperator{\Hol}{\mathsf{Hol}}
\newcommand{\lom}{\longmapsto}
\newcommand{\too}{\longrightarrow}

\newcommand{\wed}{\bigwedge}

\newcommand{\la}{\lambda}

\newcommand{\Gam}{\Gamma}

\renewcommand{\Re}{\mathbbm{Re}}
\renewcommand{\Im}{\mathbbm{Im}}

\DeclareMathOperator{\ccc}{\mathsf{c}}

\DeclareMathAlphabet{\mathscrbf}{OMS}{mdugm}{b}{n}

\DeclareMathOperator{\SO}{\mathsf{SO}}
\DeclareMathOperator{\Sp}{\mathsf{Sp}}
 \DeclareMathOperator{\SU}{\mathsf{SU}}
 \DeclareMathOperator{\GL}{\mathsf{GL}}
    
\DeclareMathOperator{\U}{\mathsf{U}}

\DeclareMathOperator{\Lie}{\mathsf{Lie}}
\DeclareMathAlphabet{\mathpzc}{OT1}{pzc}{m}{it}
\DeclareMathOperator{\Hh}{\mathsf{H}}

\DeclareMathOperator{\E}{\mathsf{E}}
\DeclareMathOperator{\Gl}{\mathsf{GL}}
\DeclareMathOperator{\Sl}{\mathsf{SL}}

\DeclareMathOperator{\Aut}{\mathsf{Aut}}
\DeclareMathOperator{\id}{\mathsf{id}}
\DeclareMathOperator{\Ss}{\mathsf{S}}
\DeclareMathOperator{\Ed}{\mathsf{End}}

\DeclareMathOperator{\J}{\mathsf{J}}

\DeclareMathOperator{\Ad}{\mathsf{Ad}}
\DeclareMathOperator{\vol}{\mathsf{vol}}
\DeclareMathOperator{\Id}{\mathsf{Id}}
\newcommand{\fr}{\mathfrak}
\newcommand{\al}{\alpha}
\newcommand{\be}{\beta}
\newcommand{\mc}{\mathcal}

\newcommand{\ep}{\varepsilon}
\newcommand{\wi}{\widetilde}
\newcommand{\cc}{\big(}

\newcommand{\rr}{\big)}

\newcommand{\om}{\omega}
\newcommand{\Om}{\Omega}

\newcommand{\lan}{\langle}
\newcommand{\ran}{\rangle}

\usepackage{amsmath}
\DeclareFontFamily{U}{mathx}{}
\DeclareFontShape{U}{mathx}{m}{n}{<-> mathx10}{}
\DeclareSymbolFont{mathx}{U}{mathx}{m}{n}
\DeclareMathAccent{\widehat}{0}{mathx}{"70}
\DeclareMathAccent{\widecheck}{0}{mathx}{"71}

\DeclareMathAlphabet{\mathscrbf}{OMS}{mdugm}{b}{n}

\newcommand{\scV}{\mathscr{V}}

\DeclareMathOperator{\ad}{ad}

\DeclareMathOperator{\Tr}{\mathsf{Tr}}

\DeclareMathOperator{\Ric}{\mathsf{Ric}}
\DeclareMathOperator{\Sca}{\mathsf{Scal}}
\DeclareMathOperator{\Hom}{\mathsf{Hom}}
\DeclareMathOperator{\ke}{\mathsf{Ker}}

\DeclareMathOperator{\dd}{d}

\newtheorem{theorem}{Theorem}[section]
\newtheorem{lem}[theorem]{Lemma}
\newtheorem{prop}[theorem]{Proposition}
\newtheorem{corol}[theorem]{Corollary}

\theoremstyle{definition}
\newtheorem{defi}[theorem]{Definition}
\newtheorem{example}[theorem]{Example}
 \newtheorem{rem}[theorem]{Remark}
 
\theoremstyle{remark}

\numberwithin{equation}{section}

\def\bex{\begin{exercise}}
\def\eex{\end{exercise}}
\def\bd{\begin{defi}}
\def\ed{\end{defi}}
\def\bt{\begin{theorem}}
\def\et{\end{theorem}}
\def\bl{\begin{lem}}
\def\el{\end{lem}}
\def\bp{\begin{prop}}
\def\ep{\end{prop}}
\def\br{\begin{rem}}
\def\er{\end{rem}}
\def\bc{\begin{corol}}
\def\ec{\end{corol}}
\def\be{\begin{example}}
\def\ee{\end{example}}
\def\pr{\begin{proof}}
\def\pro{\end{proof}}
\def\eqna{\begin{eqnarray*}}
\def\eqnaa{\begin{eqnarray}}
\def\deqna{\end{eqnarray*}}
\def\deqnaa{\end{eqnarray}}

\definecolor{dark}{rgb}{0.18,0.18,0.68}
\definecolor{mydark}{rgb}{0.78,0.08,0.08}
\definecolor{crew}{rgb}{0.2,0.5,0.2}
\definecolor{mmg}{rgb}{0.31,0.50,0.23}
\definecolor{dblue}{rgb}{0.01,0.01,0.44}
\definecolor{red}{rgb}{0.57,0.11,0.15}
\definecolor{cobalt}{RGB}{61,89,171}
\usepackage[colorlinks,citecolor=cobalt,linkcolor=cobalt,urlcolor=cobalt,pdfpagemode=UseNone,backref = page]{hyperref}

 \input ulem.sty

\title[Curvature of quaternionic skew-Hermitian manifolds]{Curvature of quaternionic skew-Hermitian manifolds and bundle constructions}

\author{Ioannis Chrysikos} 
\address{Department of Mathematics and Statistics,
Faculty of Science, Masaryk University, Kotl\' a\v rsk\' a 2, 611 37 Brno, Czech Republic}
\email{chrysikos@math.muni.cz}
\author{Vicente Cort\'es} 
\address{Department of Mathematics and Center for Mathematical Physics University of Hamburg
Bundesstraße 55, D-20146 Hamburg, Germany} 
\email{vicente.cortes@uni-hamburg.de}

\author{Jan Gregorovi\v c} 
\address{Department of Mathematics, Faculty of Science, University of Ostrava, 701 03 Ostrava, Czech Republic, and Institute of Discrete Mathematics and Geometry, TU Vienna, Wiedner Hauptstrasse 8-10/104, 1040 Vienna, Austria}
\email{jan.gregorovic@seznam.cz}

\begin{document}

\begin{abstract}
This articles is devoted to a description of the second-order differential geometry of   torsion-free almost quaternionic skew-Hermitian manifolds, 
that is, of quaternionic skew-Hermitian manifolds $(M, Q, \om)$.
We provide a curvature characterization of such integrable geometric structures,  based on the holonomy theory of symplectic connections and we
study qualitative properties of the induced Ricci tensor.
Then we proceed with bundle constructions over 
such a manifold  $(M, Q, \om)$. In particular,   we  prove the existence of almost hypercomplex skew-Hermitian structures 
on the Swann bundle over $M$ and investigate their integrability.  
 \end{abstract}

\maketitle

\tableofcontents


\section*{Introduction}\label{intro}

Hypercomplex and quaternionic geometries of \textsf{skew-Hermitian type} are geometric structures  induced by pairs $(H, \om)$, respectively $(Q, \om)$,  where $H$ is an almost hypercomplex structure,  $Q$ is an almost quaternionic structure and $\om$ is an $H$-Hermitian, respectively $Q$-Hermitian, almost symplectic 2-form.
 Such geometries  have  been  recently introduced in \cite{CGWPartI, CGWPartII} and constitute a symplectic analogue of  the better
understood (almost) hypercomplex/quaternionic Hermitians geometries.  This analogue can be encoded   in terms of $G$-structures on  $4n$-dimensional manifolds $(n>1)$, which are defined by reductions of the frame bundle to the quaternionic real form $\SO^*(2n)$ of the special  orthogonal group $\SO(2n, \C)$, and to the Lie group $\SO^*(2n)\Sp(1)$, respectively (recall   that the  structure group of
almost  hyper-Hermitian and  almost quaternionic-Hermitian structures is $\Sp(n)$ and $\Sp(n)\Sp(1)$, respectively).  The  1st-order differential geometry of such $G$-structures  was explored in
 \cite{CGWPartI, CGWPartII}, in terms of minimal connections, intrinsic torsion, and integrability conditions.

The present article is  essentially about the 2nd-order differential geometry of     quaternionic  geometries of skew-Hermitian type with a focus on the torsion-free case and on  related bundle constructions.   
 To describe our results it is convenient to assume, once and for all, that $n>1$, excluding the case when $\SO^*(2n)$ 
is compact:  $\SO^*(2)\cong\U(1)$.  We start by studying   the curvature  of torsion-free $\SO^*(2n)\Sp(1)$-structures in terms of the unique torsion-free adapted connection $\nabla^{Q, \om}$ introduced in \cite{CGWPartI}.
We do this by combining results from  \cite{CGWPartI} with  previous results on the holonomy of symplectic connections derived by   L.~Schwachh\"ofer (\cite{S1}). This allows 
the description of the curvature and Ricci tensor associated to $\nabla^{Q, \om}$ and some of their qualitative properties, see for example Corollary~\ref{QHermitianRic} and Theorem~\ref{RicciTHM}.
We also show that quaternionic skew-Hermitian symmetric spaces $G/L$ with $G$ semisimple, which are \textsf{pseudo Wolf spaces}, have necessary Hermitian Ricci tensor.   
Recall that a pseudo-Riemannian symmetric space is called a \textsf{pseudo Wolf space} if it is
of quaternionic pseudo-K\"ahler type. Such symmetric spaces were classified in \cite{AC} together with those of 
quaternionic para-K\"ahler type.
 A combination of our results 
with conclusions from \cite{AM}  allows us to derive a kind of converse of this statement, which    for the simply connected and complete case yields the classification of all torsionless $\SO^*(2n)\Sp(1)$-structures with non-degenerate $Q$-Hermitian Ricci tensor. To be more precise, we prove  that  a $4n$-dimensional quaternionic skew-Hermitian manifold $(M, Q, \om)$ with non-degenerate $Q$-Hermitian Ricci tensor  is a  (pseudo-Riemannian) quaternionic K\"ahler locally symmetric space. Therefore, when $M$ is simply connected and complete,  $(M, Q, \om)$  is one of the coset spaces  $\SU(2+p,q)/(\SU(2)\SU(p,q)\U(1))$ or $\Sl(n+1,\mathbb{H})/(\Gl(1,\mathbb{H})\Sl(n,\mathbb{H}))$. 

 The second part of this note is related to  bundle constructions associated to   a quaternionic skew-Hermitian manifold $(M^{4n}, Q, \om)$, with a focus on
the Swann bundle $\hat{M}$ over $M$ and a description of induced $\SO^*(2(n+1))$-structures. Recall  by \cite{Swann, PPS} (see also \cite{Salamon86}), that  to any  quaternionic manifold $(M, Q)$ we can associate a $\Hn^{\times}/\{\pm 1\}$-bundle over $M$ with a hypercomplex structure on the total space $\hat{M}$.
In particular, in \cite{PPS} a 1-parameter family $\hat{M}^{\ell}$ of  $\Hn^{\times}/\Z_2$-bundles was considered, each of them endowed with an almost hypercomplex structure   depending on the chosen quaternionic connection on $M$. For a specific choice of the parameter one gets a torsionless almost hypercomplex structure, i.e., a hypercomplex structure, which is actually independent 
of the initial connection. More recently in \cite{CH},  instead of a one-parameter family of bundles, it was  considered a  single principal $\Hn^{\times}/\Z_2$-bundle $\hat{M}$ over a quaternionic manifold $(M, Q)$, together
with a 1-parameter family of almost hypercomplex structures (also depending  on the initial quaternionic connection). In this case, and for the same choice of the parameter
as in \cite{PPS}, it was proved that the almost hypercomplex structure is 
$1$-integrable (and hence a hypercomplex structure) and independent of the chosen  quaternionic connection. The Swann bundle was eventually used 
by the authors of \cite{CH} in a construction of hypercomplex manifolds from quaternionic manifolds of the same dimension.

Here we mainly follow the approach of \cite{CH} and with  initial data  a quaternionic skew-Hermitian manifold $(M^{4n}, Q, \om)$, we show  that the associated Swann bundle $\hat{M}$ admits a  canonical  almost hypercomplex structure $\hat{H}$ which is in fact 1-integrable.
Then we describe associated $\SO^*(2(n+1))$-structures on $\hat{M}$, and hence we focus on almost symplectic forms on $\hat{M}$ which are $\hat{H}$-Hermitian, also called scalar 2-forms with respect to $\hat{H}$.
Here we split our study into two steps:
 First we study  the  properties of the   2-form $\hat\om=\pi^*\om$, where 
$\pi :  \hat{M}\to M$ is the bundle projection. We show that this is a basic $\hat{H}$-Hermitian presymplectic 2-form on $\hat{M}$, see Proposition~\ref{hatomega}
and Corollary~\ref{presymplectic}.  We then proceed with the classification of all the 2-forms $\beta$ on $\hat{M}$ which are $\hat{H}$-Hermitian  and vertical, which is given in Proposition \ref{32forms}. This finally  allows us to prove in Corollary~\ref{Cor_main}  the existence of canonical $\SO^*(2(n+1))$-structures on $\hat{M}$ which for non-flat manifolds $(M, Q, \om)$ are non-integrable and, more specifically, of algebraic
type $\mc{X}_{3457}$ (with respect to the classification presented in \cite[Corollary~4.17]{CGWPartI}). This algebraic type consists of those structures for which the underlying 
almost hyper-complex structure is $1$-integrable, see \cite[Theorem~1.8]{CGWPartII}. 
Hence,  in general,  our construction does not provide a  ``symplectic analogue''  of the  well-known 1-integrable $\Sp(n+1)$-structure  (i.e., a hyper-K\"ahler structure) that one obtains on $\hat{M}$ when the initial data is a 1-integrable  $\Sp(n)\Sp(1)$-structure (i.e., a  quaternionic-K\"ahler manifold).  
However, when the  initial manifold  $(M, Q, \om)$ is flat with respect to $\nabla^{Q, \om}$, i.e., $R^{Q, \om}=0$ and hence locally diffeomorphic under 
quaternionic symplectomorphisms (in terms of \cite{CGWPartI})  to the flat model $(\Hn^n\cong[\E\Hh], Q_0, \om_0)$, the situation is different. In   this case $\hat{M}$ is trivial and  we can define  $\SO^*(2(n+1))$-structures on $\hat{M}$  for which the corresponding scalar 2-form  is always symplectic, and hence of type $\mc{X}_{57}$, in general. Under 
certain natural  conditions given in Theorem \ref{theoremFlat}, we  show that  these structures are torsion-free.

\medskip
 \noindent {\bf Acknowledgments:}  V.C.\ gratefully acknowledges support by the German Science Foundation (DFG) under Germany's Excellence Strategy  --  EXC 2121 ``Quantum Universe'' -- 390833306. J.G.\  gratefully acknowledges support by Austrian Science Fund (FWF): P34369.
I. C. and J. G. thank   H. Winther (UiT Troms\o) for  useful conversations.  Some computations in this paper were performed by using 
Maple$\textsuperscript{{\tiny TM}}$. Maple is a trademark of Waterloo Maple Inc.

\section{Review of necessary theory}

\subsection{$\SO^*(2n)$-structures and  $\SO^*(2n)\Sp(1)$-structures}
 
Throughout this note, all manifolds are assumed to be smooth, connected and without boundary,
and maps are assumed to be smooth, unless otherwise mentioned. Let us fix a $4n$-dimensional  manifold  $M$.  We will denote by $\Gamma(E)$ the space of sections of a
vector bundle $\pi : E\to M$ over $M$.     

We begin with preliminaries about  $\SO^*(2n)$-structures and  $\SO^*(2n)\Sp(1)$-structures.  The  Lie group $\SO^*(2n)$ occurs as the intersection  $\Gl(n, \Hn)\cap\Sp(4n, \R)$ 
and for $n=1$  coincides with $\SO(2)=\U(1)$.  For  $n>1$  it  is  a non-compact real form of $\SO(2n, \C)$, of real dimension $n(2n-1)$,  called the \textsf{quaternionic (real) form}, while $\SO^*(2n)\Sp(1)$ is the  Lie group $\SO^{\ast}(2n)\times_{\Z_2}\Sp(1)=(\SO^{\ast}(2n)\times\Sp(1))/{\Z_2}$.    As mentioned in the introduction, we assume that $n>1$.
 The $G$-structures modeled on the above groups have been studied in \cite{CGWPartI, CGWPartII} as symplectic analogues of the more familiar almost hypercomplex/quaternionic Hermitian structures (see for example \cite{Salamon86, Swann, AM}).

 When working with this kind of $G$-structures  it is convenient to adopt the  $\E\Hh$-formalism of Salamon \cite{Salamon86}, as  in \cite{CGWPartI, CGWPartII}.    
 We consider the Lie group $\SO(2n, \C)$ of complex linear transformations preserving  
 the standard complex Euclidean metric on $\E:=\C^{2n}$ and identify $\E$ with  the 
 standard representation of $\SO^{\ast}(2n)\subset \SO(2n, \C)$. 
 Set also  $\Hh:=\C^2$ for the standard representation of $\Sp(1)$. Both  $\E, \Hh$   
 are of quaternionic type and   $\epsilon_{\E}, \epsilon_{\Hh}$ will denote    the corresponding complex anti-linear involutions.  
 Let $[\E\Hh]$ be  the real form inside $\E\otimes_\C \Hh$, fixed by the real structure $\epsilon_{\E}\otimes \epsilon_{\Hh}$. 
 The $4n$-dimensional  vector space $[\E\Hh]$ is endowed with the \textsf{standard linear quaternionic structure}, denoted by  $Q_0\cong\fr{sp}(1)$,  and defines the  standard  representation of $\SO^{\ast}(2n)\Sp(1)$,  see \cite{CGWPartI}.  
 \bd\label{defqsH}
  An   $\R$-bilinear form $h$ on $[\E\Hh]$ which is valued in endomorphisms  of $[\E\Hh]$   is called  a \textsf{quaternionic skew-Hermitian form} on $[\E\Hh]$, if
 \[
\frac12 (h(u,w)-h(w,u))\in \R\cdot \id\,, \quad \frac12 (h(u,w)+h(w,u))\in  Q_0=\sp(1)\,, 
  \] 
 and moreover $h(\J \cdot \,, \cdot)=\J\circ h(\cdot \,, \cdot)$
  for all $u, w\in [\E\Hh]$ and  $\J\in \Ss(Q_0)$,  where 
  \[
  \Ss(Q_0):=\{\mu_1 I+\mu_2 J+\mu_3 K : \mu_1^{2}+\mu_2^{2}+\mu_3^{2}=1\}\]
  is the 2-sphere of (linear) complex structures in $Q_0$.
\ed
We refer to $\Re(h)(u,w):=\frac12 (h(u,w)-h(w,u))$ (resp. $\Im(h)(u,w):=\frac12 (h(u,w)+h(w,u))$) as the real (resp.\ imaginary) part of $h$.
The   \textsf{standard  quaternionic skew-Hermitian form} on $[\E\Hh]$ is given by 
\[
h_0=\om_0\cdot\Id+\sum_{a=1}^{3}g_{0}^{a}\otimes\mc{J}_{a}
\]
with  $\Re(h_{0})=\om_0\cdot\Id$ and  $\Im(h_{0})=\sum_{a=1}^{3}g_{0}^{a}\otimes\mc{J}_{a}$, where
the triple $H_{0}:=\{\mc{J}_{a} : a=1, 2, 3 \}$ is the \textsf{standard admissible basis} of $Q_0$, $\om_0$ is a   certain non-degenerate $Q_0$-Hermitian  linear 2-form on $[\E\Hh]$  and $g_0^{a}:= \om_0(\cdot , \mc{J}_a)$.
The 2-form $\om_0$  is called the \textsf{standard scalar 2-form} on $[\E\Hh]$ and we refer to  \cite[Prop.~2.6]{CGWPartI} for its explicit expression.   
We will use the terminology \textsf{scalar 2-form}   for any non-degenerate linear 2-form  which is $Q_0$-Hermitian (and more generally, $Q$-Hermitian for some linear quaternionic structure $Q$).   
 Any other quaternionic skew-Hermitian form $h$ (respectively, scalar 2-form $\om$) on $[\E\Hh]$ is conjugated to $h_0$ (respectively, to $\om_0$), by an element in $\Gl(n,\Hn)$.
 In such terms, according to \cite[Prop.~1.15]{CGWPartI} the Lie group $\SO^*(2n)\Sp(1)$ can be viewed as the stabilizer of  $h_0$  in $\Gl(n, \Hn)\Sp(1)$.

\smallskip
Let us now discuss the manifold setting. Recall  that an \textsf{almost  hypercomplex structure} on $M$ is  a reduction of the frame bundle $\mc{F}(M)$ of $M$ to $\Gl(n, \Hn)$, that is, a $\Gl(n, \Hn)$-structure.  This means that $M$ admits  a triple $H=\{J_{a} : a=1, 2, 3\}$ of smooth endomorphisms  $J_{a}\in\Ed(TM)$ satisfying the quaternionic identity    
\[
J_{1}^2=J_{2}^2=J_{3}^2=-\Id=J_{1}J_{2}J_{3}\,.
\]
Such a pair $(M, H)$  is said to be an \textsf{almost hypercomplex manifold}.   On the other hand, an \textsf{almost quaternionic structure} on $M$  is a $G$-structure with   $G=\Gl(n,\Hn)\Sp(1)$. Geometrically,   this is given by a rank-3 smooth sub-bundle $Q\subset \Ed(TM)\cong T^*M\otimes TM$ which  locally, in a neighbourhood of each point, is generated by an almost hypercomplex structure $H=\{J_{a} : a=1, 2,3\}$.  Such a locally defined triple $H$ is called an \textsf{admissible frame (basis)} of $Q$, and  $(M, Q)$ is said to be an
\textsf{almost quaternionic manifold}.    
 Finally, an \textsf{almost symplectic structure} on  a manifold $M$  of even dimension is a $G$-structure  with $G=\Sp(\dim M,\R)\subset\Gl(\dim M, \R)$.  Thus, in this case $M$ admits a non-degenerate   2-form $\omega$, called an \textsf{almost symplectic form}, and  $(M, \omega)$ is referred to as an \textsf{almost symplectic manifold}.  Note that in our setting $\dim M =4n$ is divisible by $4$, since we will be considering manifolds admitting almost quaternionic structures.
\bd
(a) An \textsf{almost hypercomplex skew-Hermitian  structure}   on $M$ is reduction of the frame bundle $\mc{F}(M)$ of $M$ to $\SO^*(2n)$, that is, a $\SO^*(2n)$-structure. This is equivalent to  requiring that $M$ admits  a pair $(H, \omega)$, where $H=\{J_a : a=1, 2, 3\}$  is an almost hypercomplex structure on $M$ and $\omega$ is a \textsf{scalar 2-form with respect to $H$}, i.e., $\omega$  is a non-degenerate and $H$-Hermitian:
\begin{equation}\label{QHerm}
\omega(J_{a}X, J_{a}Y)=\omega(X, Y)\,,\quad a=1, 2, 3
\end{equation}
 for any $X, Y\in\Gamma(TM)$.  	We call $(M, H, \om)$  an \textsf{almost hypercomplex skew-Hermitian manifold.} \\
(b)  An \textsf{almost quaternionic skew-Hermitian    structure} on $M$ is a $\SO^*(2n)\Sp(1)$-structure. This is equivalent to  requiring that $M$  admits a pair $(Q , \omega)$, where $Q\subset\Ed(TM)$ is an almost quaternionic structure  on $M$ and $\omega$ is  a \textsf{scalar 2-form with respect to $Q$}, i.e., $\omega$  is  non-degenerate  and  $Q$-Hermitian. Hence it satisfies (\ref{QHerm})  for any  (local)  admissible base $\{J_a : a=1, 2, 3\}$ of $Q$ and   $X, Y\in\Gamma(TM)$.  
 In this case   $(M, Q, \om)$ is called an \textsf{almost quaternionic skew-Hermitian manifold.}
 \ed

\subsection{Related tensors}
Fix an almost quaternionic skew-Hermitian manifold $(M, Q, \om)$ and denote by $\pi : P\to M$  the corresponding principal $\SO^*(2n)\Sp(1)$-bundle over $M$. 
Any point in $P$ provides an identification $T_{m}M\cong[\E\Hh]$, with $m\in M$. 
The pair $(Q, \om)$ induces a volume form $\om^{2n}$ which we will denote by $\vol$. Hence $(Q, \om)$  is, in fact, a unimodular structure. 
Choosing an admissible local frame $H=\{J_{a}\}$ of $Q$ we can also define three pseudo-Riemannian metrics of signature $(2n, 2n)$,  given by 
\begin{equation}\label{metrics}
g_{J_a}(X, Y):=\omega(X,  J_{a}Y)\,,\quad a=1, 2, 3.
\end{equation}
 These metrics are in general locally defined. They are global tensors when $H$ is a global admissible frame, that is, when  $(H, \omega)$  is an    almost hypercomplex skew-Hermitian  structure on $M$.  Next we will usually  write $g_{a}$, instead of $g_{J_{a}}$. Note that  the metrics $g_{a}$ are {\it not} $H$-Hermitian, although each $g_{a}$ is $J_a$-Hermitian, i.e., $g_{a}(J_{a}X, J_{a}Y)=g_a(X, Y)$ for all $a=1, 2, 3$. 
The  normalizer of $\Sp (1)$ in the isometry group of  $g_{a}$ is  $(\SO^*(2n)\U(1))\rtimes\Z_2$, for any $a=1, 2, 3$, see \cite{CGWPartI} for details.\\ 
 Although, the individual metrics $g_a$ are frame-dependent and in general only locally defined, we have the following globally defined tensors.\\
\noindent $\bullet$  A  smooth \textsf{quaternionic skew-Hermitian form} $h$ given by $ h=\om+\sum_{a}g_{a}\otimes J_{a}$, i.e.,
 \begin{equation}\label{qsH}
h(X,Y)Z:=\omega(X,Y)Z+\sum_{a} g_{a}(X, Y)J_aZ
\end{equation}
 for any $X, Y, Z\in\Gamma(TM)$. 
This tensor is independent of the   admissible base $\{J_{a} : a=1, 2, 3\}$ of $Q$,  and the group $\SO^*(2n)\Sp(1)$ is precisely the   stabilizer of $h$ in $\Gl(n, \Hn)\Sp(1)$. 
 In particular, $h$ can be viewed as globally defined $(1, 3)$-tensor field on $M$, whose evaluation $h_{m}=\om_m+\sum_{a}(g_{a}\otimes J_{a})_{m}$ at  any $m\in M$  defines a
  quaternionic skew-Hermitian form on $T_{m}M\cong [\E\Hh]$,  in the sense of Definition~\ref{defqsH}. \\
\noindent $\bullet$  A  covariant 4-tensor $\Phi$ defined by
\[
\Phi:=\sum_{a}g_{J_{a}}\otimes g_{J_{a}}= \omega(\cdot\,, \Im(h)\cdot)\,.
\] 
This tensor  is also independent of the local admissible frame of $Q$ and in particular  $\Phi\in\Gamma(S^2S^2T^*M)$. 
We call $\Phi$ the \textsf{fundamental 4-tensor} on $(M, Q, \om)$,  since it forms the symplectic analogue of the fundamental 4-form  related to (almost) hypercomplex/quaternionic Hermitian structures. Its total symmetrization $\mathrm{Sym}(\Phi )\in \Gamma (S^4T^*M)$ defines another canonical tensor. The latter tensor was considered in \cite{CGWPartI}  (and denoted by $\Phi$ there).

To summarize, the tensors $h$ and $\Phi$   provide an alternative approach to almost quaternionic skew-Hermitian structures, and each of them  serves as a defining tensor of the $G$-structures under consideration. For instance
 \bp\label{defh}
A $4n$-dimensional connected smooth manifold $M$ $(n>1)$ admits a  $\SO^{\ast}(2n)\Sp(1)$-structure, if and only if  admits a smooth $(1, 3)$-tensor $h$ which in a local frame  of $TM$ is given by the standard quaternionic skew-Hermitian form $h_{0}$.
   \ep

\subsection{Adapted connections and 1-integrability}\label{AdaConnec}
Let $\pi : P\to M$ be a $G$-structure  on a smooth manifold $M$, where $G\subset\Gl(\scV)$ is a closed subgroup and $\scV$ is a fixed reference vector space, which we identify with $T_{m}M$ $(m\in M)$ via some coframe $u : T_{m}M\to\scV$.   An object naturally associated to $P$ is the  \textsf{first prolongation} $\fr{g}^{(1)}$ of the Lie algebra $\fr{g}=\Lie(G)$ of $G$, given by (see \cite{AM, CGWPartI} for details)
\[
\fr{g}^{(1)}:= (\scV^*\otimes\fr{g})\cap(S^2\scV^*\otimes\scV)=\{\alpha\in\scV^*\otimes\fr{g} : \alpha(u)w=\alpha(w)u, \ \forall \ u, w\in\scV\}\subset\Hom(\scV, \fr{g})\,.
\]
 $G$-structures with $\fr{g}^{(1)}=\{0\}$ are said to be of \textsf{order one}. For a $G$-structure of order one an adapted connection to the   $G$-structure, or in other words a \textsf{$G$-connection}, is entirely determined  by its torsion. Moreover,  when $\fr{g}^{(1)}=\{0\}$  a so-called  \textsf{minimal connection} of the $G$-structure (with respect to some   \textsf{normalization condition}, see \cite{AM}), must be unique, and in particular a torsion-free $G$-connection is unique. All quaternionic-like structures  (in the sense of \cite{AM}), together with $\SO^*(2n)$- and $\SO^*(2n)\Sp(1)$-structures have trivial first prolongation, except of $\Gl(n, \Hn)\Sp(1)$-structures (see  \cite{AM, CGWPartI}). It is also well-known that $G$-structures with $G\in\{\Gl(n; \C), \Sp(2n; \R)\}$ satisfy $\fr{gl}(n; \C)^{(1)}=\scV\otimes S^{2}\scV^*$ and $\fr{sp}(2n; \R)^{(1)}=S^3\scV^*$, respectively, where $\scV=\R^{2n}$ in this case.

Recall that an almost hypercomplex structure $H=\{J_{a} : a=1, 2, 3\}$ is called \textsf{hypercomplex} when it is 1-integrable, that is, there is a torsionless connection $\nabla^{H}$ with $\nabla^{H}J_{a}=0$, 
for any $a=1, 2, 3$. Then, $(M, H)$ is called  a hypercomplex manifold. It is well-known  that on a hypercomplex manifold $(M, H)$ there exists a unique such connection $\nabla^{H}$, called the \textsf{Obata connection}.  More generally, the Obata connection is defined as a canonical (but not unique) almost hypercomplex connection on any almost hypercomplex manifold. Its torsion measures the deviation from being a hypercomplex manifold. Similarly,  an almost quaternionic manifold $(M, Q)$ is said to be \textsf{quaternionic} if there exists a torsionless connection $\nabla$ which preserves $Q$, that is, $\nabla_{X}\Gam(Q)\subset\Gam(Q)$ for all $X\in\Gam(TM)$.  Such connections, i.e., \textsf{quaternionic connections} (also called \textsf{Oproiu connections}), 
subordinated to a quaternionic structure are  not unique and  form an affine space over the vector space of 1-forms on $M$. 
Note however that for unimodular quaternionic structures $(Q, \vol)$ we have  a unique adapted connection, see \cite{AM} for details.
 
Let us now  focus   on almost hypercomplex skew-Hermitian manifolds  $(M^{4n}, H, \om)$ or   almost quaternionic skew-Hermitian manifolds  $(M^{4n}, Q, \om)$, where $n>1$ for both cases. We will recall conditions encoding the 
1-integrability of such $G$-structures, as it was recently derived  in \cite{CGWPartI, CGWPartII}.  In particular, by \cite{CGWPartI}  we know explicitly the corresponding   (unique) minimal   adapted connections, and the corresponding  normalization conditions.  These adapted connections are denoted by $\nabla^{H, \omega}$ and $\nabla^{Q, \omega}$, respectively, and it will be useful to recall their  expressions.
\bt \label{connections} \textnormal{(\cite{CGWPartI})}
{\rm (1)}  The adapted connection $\nabla^{H, \om}$ on an almost hypercomplex skew-Hermitian manifold  $(M, H, \om)$ is given by $\nabla^{H, \om}=\nabla^{H}+A$,
where $\nabla^{H}$ is the (unique) Obata connection associated to $H$, and $A$ is the $(1, 2)$-tensor field on $M$ defined by  
\[
\omega\cc A(X, Y), Z\rr=\frac{1}{2}(\nabla^{H}_{X}\omega)(Y, Z)\,,\quad X, Y, Z\in\Gamma(TM)\,.
\]
\noindent{\rm (2)}  The adapted connection $\nabla^{Q, \om}$ on  an almost quaternionic skew-Hermitian manifold  $(M, Q, \om)$ is given by $\nabla^{Q, \om}=\nabla^{\vol}+A$,
where $\nabla^{\vol}\equiv\nabla^{Q, \vol}$ is the (unique)   unimodular  Oproiu  connection associated to the pair $(Q, \vol=\om^{n})$, and $A$ is the $(1, 2)$-tensor field on $M$ defined by the relation
\[
\omega\cc A(X, Y), Z\rr=\frac{1}{2}(\nabla^{\vol}_{X}\omega)(Y, Z)\,,\quad X, Y, Z\in\Gamma(TM)\,.
\]
\et

Let us finally recall the  torsion and  characterize the torsion-free case, crucial for this article.
\bp \label{torsion} \textnormal{(\cite{CGWPartII})}
\noindent {\rm (1a)} The torsion of $\nabla^{H, \omega}$ is given by
\[
T^{H, \omega}(X, Y)=T^{H}(X, Y)+\delta(A)(X, Y)\,, \quad\forall \ X, Y\in\Gamma(TM)\,,
\]
where  $\delta$  
 is the bundle map 
induced by  the Spencer operator $\delta :  \scV^{*} \otimes\fr{so}^*(2n) \to \wed^{2}\scV^{*}\otimes\scV$.  \\
\noindent {\rm (1b)}  $\nabla^{H, \omega}$ is torsion-free if and only if  
 \begin{equation}\label{omegTH}
 T^{H}=0\,,\quad \text{and}\quad  \nabla^{H}\omega=0\,.
 \end{equation}
In other words, $T^{H, \omega}=0$ if and only if $H$ is 1-integrable and $\nabla^{H, \omega}=\nabla^{H}$. Moreover, when $T^{H, \om}=0$, then  $\om$ is a symplectic form.\\
\noindent {\rm (2a)} The torsion of $\nabla^{Q, \omega}$ is given by
\[
T^{Q, \omega}(X, Y)=T^{\vol}(X, Y)+\delta(A)(X, Y)\,, \quad\forall \ X, Y\in\Gamma(TM)\,,
\]
where $\delta$ is the bundle map induced by the corresponding Spencer operator of $\fr{so}^*(2n)\oplus\fr{sp}(1)$.\\
\noindent {\rm (2b)} $\nabla^{Q, \omega}$ is torsion-free if and only if
 \begin{equation}\label{omegTHH}
 T^{Q}=0\,,\quad \text{and}\quad  \nabla^{\vol}\omega=0\,.
 \end{equation}
 In other words, $T^{Q, \omega}=0$ if and only if $Q$ is 1-integrable and $\nabla^{Q, \omega}=\nabla^{\vol}$. Moreover, when $T^{Q, \om}=0$, then  $\om$ is a symplectic  form.
 \ep
As usual, for the torsion-free    case we may omit the word ``almost'', since then our $G$-structures are 1-integrable (see \cite[p.~233]{AM}). Hence,  in this case one may speak  of \textsf{hypercomplex skew-Hermitian structures}, and \textsf{quaternionic skew-Hermitian structures}, respectively. 

\subsection{The integrability of  $\SO^*(2n)$-structures of certain algebraic type}

By \cite[Theorem 1.8]{CGWPartII} we know that on an almost hypercomplex skew-Hermitian manifold $(M, H, \om)$   the almost hypercomplex structure $H=\{J_{a} : a=1, 2, 3\}$ is 1-integrable, i.e., hypercomplex, if and only if $M$ is of type $\mc{X}_{3457}$ with respect to the algebraic types of intrinsic torsion. 
Below we present  a criterium which allows us to deduce when an almost hypercomplex skew-Hermitian manifold $(M, H, \om)$ of type $\mc{X}_{3457}$ is torsion-free.  Note that this also applies for the case where in addition we have $\dd\om=0$, i.e., for the class $\mc{X}_{57}\subset\mc{X}_{3457}$ (see \cite[Theorem 1.10]{CGWPartII}). 
\bl\label{criterium1}
Let $(M, H=\{J_{a} : a=1, 2, 3\}, \om)$ be an almost hypercomplex skew-Hermitian manifold  with $H$ hypercomplex, i.e., of type $\mc{X}_{3457}$. Then the adapted connection $\nabla^{H, \om}$ is torsion-free if and only if the Obata connection $\nabla^{H}$ is a metric connection with respect to  the induced pseudo-Riemannian metrics $g_{a}(\cdot, \cdot):=\om(\cdot, J_{a}\cdot)$, that is, $\nabla^{H}g_{a}=0$ for any $a=1, 2, 3$.
\el
\pr
We need to show that $T^{H, \om}$ vanishes identically if and only if $\nabla^{H}g_{a}=0$ for all $a\in \{ 1, 2, 3\}$.
Since $H$ is hypercomplex, the Obata connection is torsion-free,  $T^{H}=0$.  By Proposition \ref{torsion} we have  $T^{H, \om}=\delta(A)$,  which vanishes if and only if  $\nabla^{H}\om =0$. 
The latter condition is equivalent to $\nabla g_a=0$ (for one or equivalently for all $a$), since $g_a = \omega (\cdot , J_a\cdot )$ and $J_a$ is $\nabla$-parallel (and invertible).
\pro

On the other hand, by   \cite[Theorem 1.9]{CGWPartII} we know that on an almost hypercomplex skew-Hermitian manifold $(M, H, \om)$   the scalar 2-form $\om$ will be closed and hence symplectic if and only if $M$ is of type $\mc{X}_{1567}$. In this case we obtain the following criterium.
\bl\label{criterium2}
Let $(M, H=\{J_{a} : a=1, 2, 3\}, \om)$ be an almost hypercomplex skew-Hermitian manifold  with $\dd\om=0$, i.e., of type $\mc{X}_{1567}$. Then the adapted connection $\nabla^{H, \om}$ is torsion-free if and only if the endomorphism $T^{H}_{Z}:=T^{H}(Z, \cdot)$ is skew-symmetric with respect to $\om$, for any $Z\in\Gam(TM)$, and the  Obata connection $\nabla^{H}$ is a metric connection with respect to  the induced pseudo-Riemannian metrics $g_{a}(\cdot, \cdot):=\om(\cdot, J_{a}\cdot)$. 
\el
\pr
The condition $\dd\om=0$  is equivalent with  the relations (see \cite[Proposition 2.7]{CGWPartII})
\[
\fr{S}_{X, Y, Z}\om(T^{H}(X, Y), Z)=0\,,\quad\text{and}\quad \fr{S}_{X, Y, Z}(\nabla^{H}_{X}\om)(Y, Z)=0\,,
\]
for any $X, Y, Z\in\Gam(TM)$. Here $\fr{S}_{X, Y, Z}$ denotes the cyclic sum over $X, Y, Z$. If $T^{H}_{Z}$ satisfies $\om(T^{H}_{Z}X, Y)+\om(X, T^{H}_{Z}Y)=0$   (equivalently $\om(T^{H}_{Z}X, Y)+\om(T^{H}_{Y}Z,X)=0$) as an identity, then the first condition gives $T^{H}(X, Y)=0$ for any $X, Y$ and so $H$ is   hypercomplex. 
Therefore, under this assumption $M$  becomes of class $\mc{X}_{57}\subset\mc{X}_{1567}$ and we  conclude  from the previous lemma, since $\mc{X}_{57}\subset\mc{X}_{3457}$ as well.
Conversely, if $T^{H, \om}=0$ then by Proposition \ref{torsion} we have  $T^{H}=0$ and $\nabla^{H}\om=0$.  Thus $T^{H}_{Z}=0$ is trivially skew-symmetric with respect to $\om$. Moreover, as above, we see that the condition $\nabla^{H}\om=0$ implies that $\nabla^{H}g_{a}=0$.
\pro

\section{Curvature of torsion-free $\SO^*(2n)\Sp(1)$-structures}
\subsection{The   curvature tensor for the torsion-free case}
By Proposition \ref{torsion} it follows  that  the minimal connection $\nabla=\nabla^{Q, \om}$ associated to a (1-integrable) quaternionic skew-Hermitian  structure $(Q, \om)$ subordinated to a quaternionic structure $Q$ is torsion-free, and has holonomy group $\Hol(\nabla)$ contained in $\SO^*(2n)\Sp(1)$. To state a   converse of this statement we need to introduce the space of $\fr{g}$-valued formal curvature tensors, that is, the $\fr{g}$-module
\[
\mc{K}(\fr{g})=\big\{R\in\wed^{2}\scV^*\otimes\fr{g} : \fr{S}_{x, y, z}R(x, y)z=0\,, \ \text{for all} \  x, y, z\in\scV\big\}\,,
\]
where $\fr{g}=\Lie(\SO^*(2n)\Sp(1))=\fr{so}^*(2n)\oplus\fr{sp}(1)$, and $\scV=[\E\Hh]=\mathbb{R}^{4n}$. This $\fr{g}$-module fits into the exact sequence
\[
0\too\mc{K}(\fr{g})\too\wed^{2}\scV^*\otimes\fr{g}\too\wed^3\scV^*\otimes\scV\,,
\]
where the last map is given by the composition of the natural inclusion and skew-symmetrization.

\bp
Let $(M^{4n}, Q, \nabla)$ $(n>1)$ be a simply connected  quaternionic manifold with a quaternionic connection $\nabla$. Then the holonomy group $\Hol(\nabla)$ is contained in $\SO^*(2n)\Sp(1)$ 
if and only if for any $m\in M$ the curvature tensor $R^{\nabla}$ of $\nabla$ belongs to $\mc{K}(\fr{g})$. If this is the case, the connection $\nabla$ defines a 1-integrable $\SO^*(2n)\Sp(1)$-structure, subordinated to $Q$.
\ep
\pr
The proof is an application of the Ambrose-Singer theorem, compare with the proof of \cite[Thm.~2.2]{AM}.
\pro

By \cite[Main Thm.]{MS1}  we know that  the Lie group $\SO^*(2n)\Sp(1)\subset\GL(4n)$ is of special interest since it is 
an irreducible Lie subgroup which
occurs as the holonomy of a (non-Riemannian) torsion-free affine connection, which is not locally
symmetric.  In particular, $\SO^*(2n)\Sp(1)$ belongs to the list of  {\it exotic holonomies}, see \cite[p.~82]{MS1} and also to the list of {\it Segre holonomies} (\cite[p.~109]{MS1}). We summarize this important result as follows.
\bp\label{holonomy11} \textnormal{(\cite{MS1, S1})}
For $n\geq 2$ the Lie algebra of the  image of $\SO^*(2n)\Sp(1)$ via the standard representation is a (non-symmetric) irreducible Berger algebra and 
\begin{equation}\label{Holso}
\mc{K}(\fr{g})\cong\fr{g}=\fr{so}^*(2n)\oplus\fr{sp}(1)\,.
\end{equation}
\ep

We now recall  details from  \cite[Section~7]{MS1} or \cite[Section 3.2.3]{S1}, in order to specify the useful isomorphism given in (\ref{Holso}).
    Let  $\fr{g}\subset\Ed(\scV)$ be a    semisimple and irreducible subalgebra, where $\scV$ is some finite-dimensional vector space over $\R$ or $\C$. 
    Assume  that    there is a non-degenerate 2-form $\Upomega$ on $\scV$ such that $\fr{g}\subset\fr{sp}(\scV, \Upomega)$,  where $\fr{sp}(\scV, \Upomega)$ is the Lie algebra of linear symplectomorphisms of $(\scV, \Upomega)$.   By semi-simplicity, there exists a 
     non-trivial $\fr{g}$-equivariant map 
\[
\circ : S^{2}\scV^*\cong S^{2}\scV\too\fr{g}\,,\quad \circ(x\odot y)=x\circ y\,.
\]
Suppose that for any $A\in\fr{g}$ and some non-zero constant $\kappa\in\R$ the map   
\begin{equation}\label{RAxy}
R_{A} : \wed^{2}\scV\too\fr{g}\,,\quad R_{A}(x, y):=\kappa\thinspace\Upomega(x, y)A+x\circ Ay-y\circ Ax\,,\quad x, y\in\scV\,,
\end{equation}
 lies in $\mc{K}(\fr{g})$.  Then, according to \cite[Thm.~C]{MS1} or \cite[Thm.~3.6]{S1}, the assignment $A\lom R_{A}$ is an isomorphism $\fr{g}\too \mc{K}(\fr{g})$ and $\fr{g}$ is an irreducible (symplectic) Berger algebra.\footnote{Note that one may assume that $\kappa=2$ as in \cite{S1}, but it is more useful for us to  allow for an arbitrary non-zero constant $\kappa$ as in \cite[Sec.~7]{MS1}.} We will specify the data $(\Upomega, \circ)$ for $\fr{g}=\fr{so}^*(2n)\oplus\fr{sp}(1)$ below. 

\smallskip
So, let us    return back to a quaternionic skew-Hermitian manifold $(M, Q, \om)$ endowed with the torsion-free adapted connection $\nabla^{Q, \om}=\nabla^{Q, \vol}$. From now on we may set 
\[
\fr{g}:=\fr{so}^*(2n)\oplus\fr{sp}(1)\,, \quad \scV:=[\E\Hh]\cong\R^{4n}\,.
\]
According to Proposition \ref{holonomy11} the Lie algebra $\fr{g}=\fr{so}^*(2n)\oplus\fr{sp}(1)$ is an irreducible Berger algebra, and in this case we have
\[
 \Upomega=\om_{m}\quad \text{(scalar 2-form on}  \ T_{m}M\cong\scV)\,.
\]
For simplicity we will still write $\om$ instead of $\om_m$. We fix an admissible basis $\{J_a : a=1, 2, 3\}$  of  $Q\cong\fr{sp}(1)$  and consider the metrics  $g_{a}=g_{J_{a}}$.  Let 
\[
 \pi_{Z(Q)}(\om(x, -)\otimes y):=\frac{1}{4}\Big(\om(x, -)y-\sum_{a=1}^{3}g_{a}(x, -)J_{a}y\Big)\,, \quad x, y\in\scV
\]
be the projection from $\fr{gl}(\scV)$ to the centralizer $Z(Q)=\fr{gl}(n, \Hn)$ of $Q$ in $\fr{gl}(4n, \R )$ used in \cite{CGWPartI} (see also \cite[p.~214]{AM}, where $\pi_{Z(Q)}$ is denoted by $\pi_{h}$). 
Consider also the projection
\[
\pi_{Q}(\om(x, -)\otimes y):=-\frac{1}{4n}\sum_{a}\Tr(\om(x, -)\otimes J_{a}y)J_{a}
\]
 from $\fr{gl}(\scV)$ to  $Q$, see \cite[p.~214]{AM}. It is easy to see that $\pi_Q$ and  $\pi_{Z(Q)}$ are 
 unique $\GL(n,\Hn)\Sp (1)$-invariant projections to the submodules $Q, Z(Q)\subset \fr{gl}(4n, \R )$, respectively.
  In particular, they do 
 not depend on the admissible basis $\{ J_a\}$ 
 involved in the above formulas, as can be checked directly using that any two such basis are 
 conjugated by an element of $\Sp (1)$.
\bl\label{hequivlem}
Let $\circ : S^{2}\scV\too\fr{g}$ be a $\fr{g}$-equivariant map  
such that $R_A$ defined in {\rm (\ref{RAxy})} lies in $\mc{K}(\fr{g})$ for all $A\in \fr{g}$. Then  $\circ$ 
 decomposes as
\[
(x\circ y)=c_1\cdot(x\circ y)_{\fr{so}^*(2n)}+c_2\cdot(x\circ y)_{\fr{sp}(1)}\,,\quad\forall \ x, y\in\scV\,,
\]
where the $\fr{so}^*(2n)$-part $(x\circ y)_{\fr{so}^*(2n)}$ of $x\circ y$ {\rm(}respectively, the $\fr{sp}(1)$-part $(x\circ y)_{\fr{sp}(1)}${\rm)}  is given by
\eqna
(x\circ y)_{\fr{so}^*(2n)}&=&
\frac{1}{4}\Big(\om(x, -)y-\sum_{a=1}^{3}g_{a}(x, -)J_{a}y+\om(y, -)x-\sum_{a=1}^{3}g_{a}(y, -)J_ax\Big)\in\fr{so}^*(2n)\,,\\
(x\circ y)_{\fr{sp}(1)}&=&-\frac{1}{2n}\sum_{a=1}^{3}g_{a}(x, y)J_{a}\in\fr{sp}(1)\,, 
\deqna
and $c_1=2\kappa\neq 0$, $c_2=(nc_1)/2=n\kappa\neq 0$.   In particular,  the tensors on the right-hand side of these equations are independent of the admissible basis $\{J_a\}$ for $Q$.
\el
\pr
The scalar 2-form $\om$ is non-degenerate and we may identify $S^2\scV^*\cong S^2\scV$ with the Lie algebra $\fr{sp}(\scV, \om)$ by the map
\[
S^2\scV\ni x\odot y\lom f_{x\odot y}\in\fr{sp}(\scV, \om)\,,
\]
where the endomorphism $f_{x\odot y}\in\Ed(\scV)$  is defined by 
\[
f_{x\odot y}(z):=\om(x, z)y+\om(y, z)x\,,\quad x, y\in\scV\,.
\]
By Schur's lemma, the $\fr{g}$-equivariant projection $\circ : S^{2}\scV\too\fr{g}$ 
under consideration must be a sum of appropriate non-zero multiples of the projections  $\pi_{Z(Q)}$ and $\pi_{Q}$ 
restricted to $\fr{sp}(\scV, \om)\subset\fr{gl}(\scV)$.    In particular, the $\fr{so}^*(2n)$-part  
of the equivariant projection $\circ : S^{2}\scV\too\fr{g}$ appears by applying  $\pi_{Z(Q)}$ to the endomorphism  $f_{x\odot y}=\om(x, -)\otimes y+\om(y, -)\otimes x$, 
\[
x\odot y\lom 
(x\circ y)_{\fr{so}^*(2n)}:=\pi_{Z(Q)}(f_{x\odot y})\,,\quad x, y\in\scV\,.
\]
This   gives the first stated expression.
Similarly, to obtain the $\fr{sp}(1)$-part of  the equivariant projection $\circ : S^{2}\scV\too\fr{g}$   
we apply $\pi_{Q}$  to the endomorphism  $f_{x\odot y}=\om(x, -)\otimes y+\om(y, -)\otimes x$:
\[
x\odot y\lom 
(x\circ y)_{\fr{sp}(1)}:=\pi_{Q}(f_{x\odot y})\,,\quad x, y\in\scV\,.
\]
Having in mind the relation $\Tr(\om(x, -)\otimes y)=\omega(x,y)$, this gives
\eqna
\pi_{Q}(f_{x\odot y})&=&\pi_{Q}(\om(x, -)\otimes y)+\pi_{Q}(\om(y, -)\otimes x)\\
&=&-\frac{1}{4n}\sum_{a}\Tr(\om(x, -)\otimes J_{a}y)J_{a}-\frac{1}{4n}\sum_{a}\Tr(\om(y, -)\otimes J_{a}x)J_{a}\\
&=&-\frac{1}{4n}\sum_{a}\om(x,  J_{a}y)J_{a}-\frac{1}{4n}\sum_{a}\om(y,  J_{a}x)J_{a}\\
&=&-\frac{1}{4n}\sum_{a}g_{a}(x, y)J_{a}-\frac{1}{4n}\sum_{a}g_{a}(y, x)J_{a}\\
&=&-\frac{1}{2n}\sum_{a}g_{a}(x, y)J_{a}\,.
\deqna
 The fact that the projections $\pi_Q$ and $\pi_{Z(Q)}$ are independent of the chosen admissible basis $\{ J_a\}$  
implies that $(x\circ y)_{\fr{so}^*(2n)}$ and $(x\circ y)_{\fr{sp}(1)}$  are independent as well.

We are now ready to specify  $c_1, c_2$. To find them we rely on the first Bianchi identity.  For the map $R_A=R_{A}^{Q, \om} : \wed^{2}\scV\too\fr{g}$ compare equation~\eqref{RAxy}, we have 
 \begin{equation}\label{RQomA1}
R^{Q, \om}_{A}(x, y)=\kappa\cdot\om(x, y)A+c_{1}\cdot(x\circ Ay-y\circ Ax)_{\fr{so}^*(2n)}+c_{2}\cdot(x\circ Ay-y\circ Ax)_{\fr{sp}(1)}\,.
\end{equation}
Hence for the associated $(1, 3)$-tensor one obtains
 \eqna
R^{Q, \om}_{A}(x, y)z&=&\kappa\thinspace\om(x, y)Az+\frac{c_1}{4}\Big(\om(x, z)Ay-\sum_{a=1}^{3}g_{a}(x, z)J_{a}Ay+\om(Ay, z)x-\sum_{a=1}^{3}g_{a}(Ay, z)J_ax\Big)\\
&&-\frac{c_1}{4}\Big(\om(y, z)Ax-\sum_{a=1}^{3}g_{a}(y, z)J_{a}Ax+\om(Ax, z)y-\sum_{a=1}^{3}g_{a}(Ax, z)J_ay\Big)\\
&&-\frac{c_2}{2n}\sum_{a=1}^{3}\big(g_{a}(x, Ay)-g_{a}(y, Ax)\big)J_{a}z\,,
\deqna
for any $x, y, z\in\scV$, $A\in\fr{g}$ and some non-zero $\kappa, c_1, c_2\in\R$. Now, a short computation relative to the first Bianchi identity shows the following. 
\begin{eqnarray}
0&=&R^{Q, \om}_{A}(x, y)z+R^{Q, \om}_{A}(y, z)x+R^{Q, \om}_{A}(z, x)y\nonumber\\
&=&\fr{S}_{x, y, z}\Big(\kappa\thinspace\om(x, y)Az+\frac{c_1}{4}\om(y, x)Az-\frac{c_1}{4}\om(x, y)Az\Big)\nonumber\\
&&+\fr{S}_{x, y, z}\Big(-\frac{c_1}{4}\sum_{a}g_{a}(Ay, z)J_{a}x-\frac{c_2}{2n}\sum_{a}g_{a}(y, Az)J_{a}x+\frac{c_2}{2n}\sum_{a}g_{a}(z, Ay)J_{a}x\nonumber\\
&&+\frac{c_1}{4}\sum_{a}g_{a}(Az, y)J_{a}x\Big)\nonumber\\
&=&\fr{S}_{x, y, z}\Big(\big(\kappa-\frac{c_1}{2}\big)\om(x, y)Az+ \big(\frac{c_1}{4}-\frac{c_2}{2n}\big)\sum_{a}\big(g_{a}(Ay, x)-g_{a}(Ax, y)\big)J_{a}z\Big)\,,\label{1stbianchi}
\end{eqnarray}
for any $x, y, z\in\scV$ and $A\in\fr{g}$, where as before $\fr{S}_{x, y, z}$ denotes the cyclic sum over $x, y, z$.  Thus obviously, for the given values of $c_1, c_2$ in the statement the first Bianchi identity is satisfied. Hence these conditions are sufficient. To prove that are also necessary observe that  we can take $x,y$ in the $\omega$-orthogonal quaternionic complement of the line spanned by $z$, since 
the (real) dimension is at least $8$. For such a triple of vectors $(x,y,z)$ the only non-zero terms  in \eqref{1stbianchi}
are multiples of  $Az$ and $J_a z$, $(a=1,2,3)$.  Since the equation \eqref{1stbianchi} holds for all $A\in\fr{so}^*(2n)$, it is easy to see
that the coefficients $\kappa -\frac{c_1}{2}$ and $\frac{c_1}{4} -\frac{c_2}{2n}$ have to vanish.
\pro

Two important corollaries of Lemma \ref{hequivlem}  read as follows.
\bc\label{curvatureRQom}
 For any $x, y, z\in\scV=[\E\Hh]\cong T_{m}M$, $A\in\fr{so}^*(2n)\oplus\fr{sp}(1)$ and some non-zero $\kappa\in\R$ we have 
\eqnaa
R^{Q, \om}_{A}(x, y)z&=&\kappa\thinspace\om(x, y)Az+\frac{\kappa}{2}\Big(\om(x, z)Ay-\sum_{a=1}^{3}g_{a}(x, z)J_{a}Ay+\om(Ay, z)x-\sum_{a=1}^{3}g_{a}(Ay, z)J_ax\Big)\nonumber\\
&&-\frac{\kappa}{2}\Big(\om(y, z)Ax-\sum_{a=1}^{3}g_{a}(y, z)J_{a}Ax+\om(Ax, z)y-\sum_{a=1}^{3}g_{a}(Ax, z)J_ay\Big)\nonumber\\
&&-\frac{\kappa}{2}\sum_{a=1}^{3}\big(g_{a}(x, Ay)-g_{a}(y, Ax)\big)J_{a}z\,.\label{RQom}
\deqnaa
\ec

\bc\label{adjointbundle}
Set $G=\SO^*(2n)\Sp(1)$, and let $(M,  Q, \om, \nabla^{Q, \om})$ be a  quaternionic skew-Hermitian manifold with corresponding principal $G$-bundle $\pi : P\to M$. Then, there is a section $A$ of the adjoint bunde $\fr{g}_{P}=P\times_{G}\fr{g}\subset T^*M\otimes TM$   such that curvature  $R^{Q, \om}$ of $\nabla^{Q, \om}$ (which can be seen as a section of $\wed^{2}T^*M\otimes\fr{g}_{P}$) corresponds to  the section $R^{Q, \om}_{A}$ of  the bundle over $M$ with fiber $\mc{K}(\fr{g})$.   
\ec
\begin{proof}This follows from the isomorphism $\mathfrak{g} \to \mc{K}(\fr{g}), A \mapsto R_A$. 
\end{proof}

\subsection{The   Ricci tensor for the torsion-free case}
Let us now pass to the 
  Ricci tensor of $R^{Q, \om}_{A}$, which  is  traditionally defined by   (see for example \cite[p.~222]{AM})
\[
\Ric^{Q, \om}_{A}(y, z):=\Tr\big\{x\lom R^{Q, \om}_{A}(x, y)z\big\}\,,\quad x, y, z\in\scV\,.
\]
\bp\label{RicQomProp}
For any $y, z\in\scV=[\E\Hh] \cong T_{m}M$ and $A\in\fr{so}^*(2n)\oplus\fr{sp}(1)$ the Ricci tensor associated to $R^{Q, \om}_{A}$ is given by
\eqnaa
\Ric^{Q, \om}_{A}(y, z)&=&(2n+1)\kappa\thinspace\om(Ay, z)+  \frac{\kappa}{2}\sum_{a}g_{a}(y, z)\Tr(J_{a}A) 
-\kappa\sum_{a}\om(J_{a}AJ_{a}y, z)\label{RicAyzsimple}\,.
\deqnaa
\ep
\pr
Since $J_{a}\in\fr{sp}(4n, \R)$ we have $\Tr(J_a)=0$. Similarly, any $A\in\fr{g}$  satisfies $\Tr(A)=0$, since $\fr{g}\subset\fr{sp}(4n, \R)$.  By Corollary \ref{curvatureRQom} this implies that  the traces of the 5th and 6th term  in (\ref{RQom})  vanish. Therefore, we obtain
\eqna
\Ric^{Q, \om}_{A}(y, z)&=&\kappa\thinspace\om(Az, y)+\frac{\kappa}{2}\big(\om(Ay, z)-\sum_{a}g_{a}(J_{a}Ay, z)+4n\om(Ay, z)\big)\\
&&-\frac{\kappa}{2}\big(\om(Ay, z)-\sum_{a}g_{a}(y, z)\Tr(J_{a}A)-\sum_{a}g_{a}(AJ_{a}y, z)\big)\\
&&-\frac{\kappa}{2}\sum_{a}\big(g_{a}(J_{a}z, Ay)-g_{a}(y, AJ_{a}z)\big)\,,\\ 
\deqna
for any  $y, z\in\scV$. 
In this relation the 2nd term and the 5th  term cancel  each other.  Note also that its 6th term is {\it not} zero, in general.  In fact,  
 below we will show that for $A=\sum_{b=1}^{3}c_{b}J_{b}\in\fr{sp}(1)\setminus \{ 0\}$ the bilinear form $(y,z) \mapsto \sum_{a}g_{a}(y, z)\Tr(J_{a}A)$ is not zero.
 Moreover,     since $\om$ is $\fr{g}$-invariant, 
the sum of its 1st and  4th term  gives $(2n+1)\kappa\thinspace\om(Ay, z)$.  Thanks to these conclusions we deduce that  
\eqna
\Ric^{Q, \om}_{A}(y, z)&=&(2n+1)\kappa\thinspace\om(Ay, z) -\frac{\kappa}{2}\sum_{a}g_{a}((J_{a}A-AJ_{a})y, z)+\frac{\kappa}{2}\sum_{a}g_{a}(y, z)\Tr(J_{a}A)\nonumber\\
&&-\frac{\kappa}{2}\sum_{a}\big(g_{a}(J_{a}z, Ay)-g_{a}(y, AJ_{a}z)\big)\,. 
\deqna
 To arrive  at the given expression (\ref{RicAyzsimple}) we use the definition  $g_{a}(u, w)=\om(u, J_{a}w)$, $u, w\in\scV$, 
 and  the fact that $\om$ is $Q$-Hermitian  (see \cite[Prop.~1.10]{CGWPartI}),  that is,  $J_a^*\om = \om$ (or equivalently $J_{a}\in\fr{sp}(4n, \R)$) for all $a\in \{ 1,2,3\}$.

Then,  setting ${\sf s}:= -\frac{\kappa}{2}\sum_{a}g_{a}((J_{a}A-AJ_{a})y, z) -\frac{\kappa}{2}\sum_{a}\big(g_{a}(J_{a}z, Ay)-g_{a}(y, AJ_{a}z)\big)$ 
we compute 
\eqna
{\sf s}&=&\frac{\kappa}{2}\sum_{a}\big(\om(AJ_{a}y, J_{a}z)+\om(y, J_{a}AJ_{a}z)\big)= -\frac{\kappa}{2}\sum_{a}\big(\om(J_{a}AJ_{a}y, z)+\om(J_{a}y, AJ_{a}z)\big)\\
&=& -\frac{\kappa}{2}\sum_{a}\big(\om(J_{a}AJ_{a}y, z)+\om(J_{a}AJ_{a}y, z)\big)=-\kappa\sum_{a}\om(J_{a}AJ_{a}y, z)\,,
\deqna
and this proves (\ref{RicAyzsimple}). 
\pro
\bc\label{QHermitianRic}
$(1)$ For $A\in\so^*(2n)$ the Ricci tensor satisfies 
\[
\Ric^{Q, \om}_{A}(y, z)=2(n+2)\kappa\thinspace\om(Ay, z)\,,\quad y, z\in\scV\,.
\]
$(2)$ For $A\in\sp(1)$ the Ricci tensor satisfies 
\[
\Ric^{Q, \om}_{A}(y, z)=4n\,\kappa\thinspace\om(Ay, z)\,,\quad y, z\in\scV\,.
\]
$(3)$ The Ricci tensor $\Ric^{Q, \om}_{A}$ is $Q$-Hermitian if and only if $A\in\fr{so}^*(2n)$.\\
$(4)$ The Ricci tensor $\Ric^{Q, \om}_{A}$ is symmetric for any $A\in\fr{so}^*(2n)\oplus\fr{sp}(1)$.
\ec
\pr
To prove $(1)$, we first observe that the second term in (\ref{RicAyzsimple}) vanishes, since 
$\Tr (J_{a}A)=0$  for all $A\in\fr{so}^*(2n)$ (see \cite{CGWPartI}). For   $A\in\fr{so}^*(2n)$ we also 
have that $AJ_{a}=J_{a}A$ and, thus, $J_aAJ_a = J_a^2A=-A$.
Therefore,  the last sum in (\ref{RicAyzsimple}) gives
\[
-\kappa\sum_{a}\om(J_{a}AJ_{a}y, z)=3\kappa\,\om(Ay, z)\,,
\]
and hence,   $\Ric^{Q, \om}_{A}(y, z)=(2n+1)\kappa\thinspace\om(Ay, z)+3\kappa\thinspace\om(Ay, z)=(2n+4)\kappa\thinspace\om(Ay, z)\,.$

\noindent To prove $(2)$, we can assume without loss of generality that $A=J_1$. In fact, we can always choose the admissible basis $(J_a)$ such that 
$A$ is a multiple of $J_1$, and then we may as well assume that $A=J_1$. Thus we obtain
\begin{eqnarray*}
 (2n+1)\kappa\thinspace\om(Ay, z) &=& -(2n+1)\kappa\, g_1(y,z)\,,\\
\frac{\kappa}{2}\sum_{a}g_{a}(y, z)\Tr(J_{a}A) &=& \frac{\kappa}{2}g_{1}(y, z)(-4n)=-2n\,\kappa\, g_{1}(y, z)\,,
\end{eqnarray*}
and 
\[
 -\kappa\sum_{a}\om(J_{a}AJ_{a}y, z) = -\kappa \big(\om(J_{1}^3y, z)- \omega (J_2^2J_1y,z) -\omega (J_3^2J_1y,z)\big)=-\kappa\, \omega (J_1y,z)=\kappa\, g_1(y,z)\,.
\]
Finally, these expressions sum up to $\Ric^{Q, \om}_{A}(y, z)= \Ric^{Q, \om}_{J_1}(y, z)= -4n\,\kappa\,g_1(y,z)= 4n\, \omega (Ay,z)$, which proves (2). 
Finally, for $(3)$ we  recall that $\Ric_{A}^{Q, \om}$ is   $Q$-Hermitian when
\[
\Ric^{Q, \om}_{A}({\sf J}x, {\sf J}y)=\Ric^{Q, \om}_{A}(x, y)\,,
\]
for all $x, y\in\scV$, ${\sf J}\in S^{2}(Q)$.  This is equivalent to saying that  $\Ric^{Q, \om}_A$
 is $\fr{sp}(1)$-invariant and the conclusion follows  from the expressions for $\Ric^{Q, \om}_{A}$ given in the parts  $(1)$ and $(2)$, when $A\in\fr{so}^*(2n)$ and $A\in\fr{sp}(1)$, respectively. Indeed, according to (\ref{RicAyzsimple}) the Ricci tensor depends linearly on $A$ and since $\fr{g}=\fr{so}^*(2n)\oplus\fr{sp}(1)$, by $(1)$ and $(2)$ we can write
\[
\Ric^{Q, \om}_{A}(x, y)=(2n+4)\kappa\thinspace\om(A_1x, y)+4n\,\kappa\,\om(A_2x, y)\,, 
\]
with $A_1\in\fr{so}^*(2n)$ and $A_2\in\fr{sp}(1)$ such that  $A=A_1+A_2$.  Thus $\Ric^{Q, \om}_{A}$ is $\sp(1)$-invariant if and only if $A=A_1\in\fr{so}^*(2n)$. The final claim about the symmetry  of $\Ric^{Q, \om}_{A}$ is obvious from 
parts (1) and (2), in view of the fact that $A$ is $\omega$-skew-symmetric.
\pro

Now, for the manifold setting we  adopt the following definition.
 \bd
Let $(M, Q, \om, \nabla^{Q, \om})$ be a quaternionic skew-Hermitian manifold. The Ricci tensor $\Ric^{Q, \om}=\Ric^{\nabla^{Q, \om}}$  associated to $\nabla^{Q, \om}$ is called  \textsf{$Q$-Hermitian} if  at any point $m\in M$ it satisfies the relation
\[
\Ric^{Q, \om}_{m}({\sf J}x, {\sf J}y)=\Ric_{m}^{Q, \om}(x, y)
\]
for any $x, y\in\scV\cong T_{m}M$ and any ${\sf J}\in S^{2}(Q)$.
\ed
According to \cite[Thm.~5.3]{AM} the Ricci tensor $\Ric^{\nabla}$ of a quaternionic connection $\nabla$ on a quaternionic manifold $(M, Q)$ is $Q$-Hermitian if and only if  the curvature tensor $R^{\nabla}$ is $\Sp(1)$-invariant, or equivalently, $\fr{sp}(1)$-invariant.
 Note that $R^{\nabla}$ is $\sp(1)$-invariant if and only if  
 \[
(b\thinspace R^{\nabla}_{m})(x, y, z, w)=0, 
 \]
 for all $m\in M$, $x, y, z, w\in\scV\cong T_{m}M$ and $b\in Q_m$,
 where the Lie algebra $Q_m$ acts by derivations on the tensor algebra, such that
 \[
 b\thinspace R^{\nabla}_{m}(x, y, z, w)=-R^{\nabla}_{m}(bx, y, z, w)-R^{\nabla}_{m}(x, by, z, w)-R^{\nabla}_{m}(x, y, bz, w)-R^{\nabla}_{m}(x, y, z, bw)\,.
 \]
 Since $\nabla^{Q, \om}$ is a quaternionic connection (which is also symplectic), by Proposition \ref{adjointbundle} and Corollary~\ref{QHermitianRic} we deduce that
 
\bt\label{RicciTHM} Let $(M, Q, \om)$ be a quaternionic skew-Hermitian manifold   endowed with the  torsion-free $\SO^*(2n)\Sp(1)$-connection  $\nabla=\nabla^{Q, \om}=\nabla^{Q, \vol}$ and let $A$ be the unique section of the adjoint bundle $\mathfrak{g}_P$ for which $R_{A}^{Q, \om}=R^{Q, \om}=R^\nabla$, where $P\rightarrow M$ denotes the $\SO^*(2n)\Sp(1)$-structure.   Then the following are equivalent:\\
$(1)$ $A\in\fr{so}^*(2n)$.\\
$(2)$ $(b\thinspace R_{A}^{Q, \om})(X, Y, Z, W)=0$, for any $X, Y, Z, W\in\Gam(TM)$  and $b\in\Gam(Q)$.\\
$(3)$ The Ricci tensor associated to $\nabla^{Q, \om}$ is $Q$-Hermitian.\\
$(4)$ The Ricci tensor associated  to $\nabla^{Q, \om}$ can be written as $\Ric^{Q, \om}(X, Y)=2(n+2)\kappa\thinspace\om(AX, Y)$, 
for any $X, Y\in\Gam(TM)$, where $\kappa$ is the non-zero constant appearing in \eqref{RQom}. 
 \et

 \subsection{Torsionless $\SO^*(2n)\Sp(1)$-structures with non-degenerate $Q$-Hermitian Ricci tensor}
 A natural source where one may  look for  examples of $\SO^*(2n)\Sp(1)$-structures admitting a 
 torsion-free compatible connection 
   is the category of symmetric spaces $G/L$ with $G$ semisimple.  
   Such spaces were classified in  \cite{CGWPartI}. Up to covering,  they are exhausted  by the   following families of
  symmetric spaces:
\[
\SO^*(2n+2)/\SO^*(2n)\U(1)\,,\; \SU(2+p,q)/(\SU(2)\SU(p,q)\U(1))\,,\;\ \Sl(n+1,\mathbb{H})/(\Gl(1,\mathbb{H})\Sl(n,\mathbb{H}))\,.
\]
While the last two are  pseudo Wolf spaces (see \cite{AC}), the first symmetric space   is not, as explained in Example~\ref{notWolf} (see also \cite[p.~2654]{CGWPartI}).
At  present, these three families  constitute the only known global examples of quaternionic skew-Hermitian manifolds. 
 In particular, the  construction of  global examples of non-symmetric torsion-free $\SO^*(2n)\Sp(1)$-structures is an open problem.  Below we analyze some of their curvature features, in terms of the previous paragraph. As we will see, 
this will be helpful to obtain  a classification result about torsionless $\SO^*(2n)\Sp(1)$-structures with non-degenerate $Q$-Hermitian Ricci tensor.

\be\label{notWolf}
Let us first consider the  space  $M=G/L=\SO^*(2n+2)/\SO^*(2n)\U(1)$. 
 This 
 is a member of the  2-parameter family of symmetric spaces $\SO^*(2n+2m)/\SO^*(2n)\SO^*(2m)$ (it occurs for $m=1$, see \cite[p.~414]{Gil}).  According to  \cite[Thm.~3.9]{Pont} or \cite[Example 4.4.1]{Jan},   $M$ admits an invariant quaternionic structure $Q$. 
By \cite{Pont} we know that    $M$ carries an invariant complex structure $J$ which comes from a global section $M\to Z$, where $Z$ is the twistor space over $M$.
Let $\fr{g}=\fr{l}\oplus\fr{m}$ be the Cartan decomposition associated to the Lie algebra 
  $\fr{g}=\fr{so}^*(2n+2)$  of $G=\SO^*(2n+2)$, where $\fr{l}=\Lie(L)$ is the Lie algebra of $L$  and let us identify   $\fr{m}\cong T_{eL}(G/L)$.  The   isotropy representation $\chi : L\to\Aut(\fr{m})$ is the restriction of the standard representation $[\E\Hh]$ of $\SO^*(2n)\Sp(1)$  to $\SO^*(2n)\U(1)$, which is irreducible as a consequence of the irreducibility of $[\E\Hh]$ as an $\SO^*(2n)$-module.
    The latter follows since  the irreducible complex 
 $\SO^*(2n)$-module $\E$ is of quaternionic type and therefore irreducible also as a real module and moreover we have $[\E\Hh]\cong \E$, as 
 $\SO^*(2n)$-modules. 
Since $M$ is isotropy irreducible,  we  can identify (up to a scaling) the complex structure $J$   with the invariant complex structure induced by the $\Ad(L)$-invariant  endomorphism $J_{o} : \fr{m}\to\fr{m}$,  defined by   $ J_{o}:=\chi_{*}(\Upupsilon)=\ad(\Upupsilon)$ 
   for some $\Upupsilon\in\fr{u}(1)$.
    Then,  the   $\SO^*(2n+2)$-invariant quaternionic skew-Hermitian structure induced by the embedding  $\SO^*(2n)\U(1)\subset \SO^*(2n)\Sp(1)$ can be described in terms of the 
$\Ad(L)$-invariant  (symmetric) pseudo-Hermitian metric $g_{o}=\langle \ , \ \rangle_{\fr{m}}$  on $\fr{m}$ corresponding to $J_{o}$, as follows. Set     
 \[
 \omega_o(X, Y):=\langle J_{o}X, Y\rangle_{\fr{m}}\,, \quad X, Y\in\fr{m}\,.
 \]
It is easy to see that $\om_o$ is  an $\Ad(L)$-invariant  non-degenerate  2-form on $\fr{m}$ which is $Q_{o}$-Hermitian, that is,  an invariant linear scalar 2-form on $(\fr{m}, Q_{o})$. Hence the corresponding $\SO^*(2n+2)$-invariant 2-form $\om$ on $M=G/L$  is $Q$-Hermitian, and by Schur's lemma this is the unique (up to scale) $\SO^*(2n+2)$-invariant scalar 2-form on $(M=G/L, Q)$.  Consequently, $(M, Q, \om)$ is a  homogeneous almost quaternionic skew-Hermitian manifold.  We have $g_{J_{o}}=g_{o}=\langle \ , \ \rangle_{\fr{m}}$ and this   induces the unique   invariant Einstein metric on $M=G/L$ of signature $(2n, 2n)$, which is a multiple of the Killing form of $\SO^*(2n+2)$ restricted to $\fr{m}$.  
Since the pair $(Q, \om)$ lives  in a symmetric space,  it is  actually  a quaternionic skew-Hermitian structure (torsion-free). This is because   the adapted connection  $\nabla^{Q,\omega}$  coincides with the canonical connection $\nabla^{0}$ on  $M=G/L$, and the latter is torsion-free since $[\fr{m}, \fr{m}]\subset\fr{l}$.    Moreover, since $\nabla^0$ is torsion-free and unimodular, the corresponding Ricci tensor  $\Ric^{0}$ is symmetric. We will show however that  $\Ric^0$ is {\it not} $Q$-Hermitian.  
 We extend $J_{o}$ to an admissible basis $\{I, J_{o}, IJ_{o}\}$ of $Q_o$.  Since $\lan X, Y\ran_{\fr{m}}=\om_{o}(X, J_{o}Y)$ for any $X, Y\in\fr{m}$ we have
\eqna
\Ric^{0}(IX, IY)&=&c\thinspace\om_o(IX, J_{o}(IY))=-c\thinspace\om_{o}(IX, IJ_{o}Y)=-c\thinspace\om_{o}(X, J_{o}Y)=-\Ric^{0}(X, Y)\,,\\
\Ric^{0}(J_{o}X, J_{o}Y)&=&c\thinspace\om_{o}(J_oX,  J^2_{o}Y)=c\thinspace\om_{o}(X, J_{o}Y)=\Ric^{0}(X, Y)\,,
\deqna
for all $X, Y\in\fr{m}$, and  some constant $c$ (Einstein constant), where we have used the property that $\om_o$ is $Q_o$-Hermitian.
So, we see that $J_o$ is the only complex structure in $Q_o$ which preserves $\Ric^0|_o=cg_o$, $o=eL$. 
In particular, $\Ric^0$ is {\it not} $Q$-Hermitian (and the pseudo-Riemannian symmetric $M=G/L$ is not of quaternionic 
pseudo-K\"ahler type). By Corollary \ref{adjointbundle} and Corollary \ref{QHermitianRic}  we have $\Ric^{0}(X, Y)=4n\,\kappa\,\om_{o}(AX, Y)$,  for some $A\in\fr{sp}(1)$ and  real number $\kappa\neq 0$. Indeed we see that $A$ is a multiple of $J_{o}\in\fr{u}(1)$, 
\[
A=-\frac{c}{4n\kappa}J_{o}\,.
\]
Note that  the Einstein constant $c=\frac{\Sca^{g_{o}}}{4n}$ is negative if we use the Killing form 
(or a positive multiple) thereof to define the metric.
   \ee
   \br   By the above description it follows that  $M=G/L$ is an example of  quaternionic manifold with a 
   quaternionic connection of holonomy $\SO^*(2n)\U(1)$. The latter manifolds do always admit a complex structure subordinate to the quaternionic structure (corresponding to $J_0$ in the example) as can be read off from their holonomy group. As observed in 
   \cite[p.\ 176]{Pont},  
there is no parallel hypercomplex structure with respect to the given quaternionic connection.
\er
 Independently of the particular 
quaternionic connection we can show the following proposition.
\bp Let $(M,Q)$ be a quaternionic manifold that admits a quaternionic connection $\nabla$ with holonomy algebra in $\so^*(2n)\oplus \mathfrak{u}(1)$  
such that the projection to  $\mathfrak{u}(1)$ is non trivial. Then, there is no local hypercomplex structure subordinate to $Q$ which extends the parallel complex structure. 
\ep 
\begin{proof} From the assumption on the holonomy we obtain the existence of a local admissible frame $\{J_a\}$ of $Q$ 
extending the parallel complex structure $J_1$ 
and three local one-forms $(\omega_a)$ 
such that 
\[
 \nabla J_a = \omega_b \otimes J_c -\omega_c\otimes J_b\,,
\] 
for any cyclic permutation of $\{ 1,2,3\}$ and $\omega_2=\omega_3=0$. 
It follows from \cite[Prop.\ 3.2]{AM} that the almost hypercomplex structure $\{J_a : a=1, 2, 3\}$ is integrable if and 
only if $\omega_a \circ J_a$ is independent of $a\in \{ 1,2,3\}$. This is clearly not the case. 
\end{proof}

\be\label{qsHpWolf}
We now concentrate on the pseudo Wolf spaces  (see \cite{AC})
\[
\SU(2+p,q)/(\SU(2)\SU(p,q)\U(1))\,,\quad\ \Sl(n+1,\mathbb{H})/(\Gl(1,\mathbb{H})\Sl(n,\mathbb{H}))\,.
\]
 These symmetric spaces $G/L$ satisfy $\chi(L)\cap\Sp(1)=\Sp(1)$, where $\chi : L\to\Aut(\fr{m})$ is the corresponding isotropy representation, and this induces a quaternionic structure $Q$ which is $G$-invariant.   In this case the Ricci tensor corresponding to $\nabla^{Q, \om}$ turns out to be $Q$-Hermitian. We describe the details for   the homogeneous space $M=G/L=\SU(2+p,q)/(\SU(2)\SU(p,q)\U(1))$ and  similarly is treated the other coset. The isotropy action of $\U(1)$ gives rise to an invariant complex structure $I$ on $\SU(2+p,q)/(\SU(2)\SU(p,q)\U(1))$, satisfying $I\notin\Gam(Q)$. Let $g_{o}=\lan \ , \ \ran_{\fr{m}}$ be the $\SU(2+p,q)$-invariant pseudo quaternionic K\"ahler metric on $M$, which is a multiple of the Killing form of $\SU(2+p,q)$ restricted to $\fr{m}$, where   $\fr{m}\cong T_{o}G/L$ is the reductive complement associated to the symmetric reductive decomposition of $\fr{su}(2+p, q)$. Recall that $g$ is a $G$-invariant Einstein metric on $M$, unique up to scale.
Set 
\[
\om_{o}(X, Y):=\lan IX, Y\ran_{\fr{m}}\,, \quad X, Y\in\fr{m}\,.
\]
This defines a non-degenerate 2-form on $\fr{m}$ which is $\Ad(L)$-invariant  and $Q_{o}$-Hermitian. Hence it induces a $G$-invariant scalar 2-form $\om$ on $M=G/L$, which by Schur's lemma is unique (up to scale). Thus   $(Q, \om)$ is a  $G$-invariant quaternionic skew-Hermitian structure on $M$.    As in the previous example, we have $\nabla^{Q, \om}=\nabla^{0}$.  Let $\{J_a : a=1, 2, 3\}$ be an admissible basis of $Q$.  Suppose that 
  $IJ_{a}=-J_{a}I$ for some $a$. Then, the  triple $\{I, J_{a}, IJ_{a}\}$ will be another admissible basis of $Q$, a contradiction since $I\notin\Gamma(Q)$. Hence $IJ_{a}=J_{a}I$ for any $a=1, 2, 3$ and since $\om$ is $Q$-Hermitian we deduce  that
  \eqna
    \Ric^{0}(J_{a}X, J_{a}Y)&=&c\thinspace\lan J_{a}X, J_{a}Y\ran_{\fr{m}}=-c\thinspace\om(IJ_{a}X, J_{a}Y)=-c\thinspace\om(J_{a}IX, J_{a}Y)\\
    &=&-c\thinspace\om(IX, Y)=c\thinspace\lan X, Y \ran_{\fr{m}}=\Ric^{0}(X, Y)\,,
  \deqna
 for any $X, Y\in\fr{m}$. It follows that $\Ric^{0}$ is $Q$-Hermitian.  In particular, the matrix $A\in\fr{so}^*(2n)$ given in Theorem  \ref{RicciTHM} is  a multiple of the complex structure $I$, 
 \[
 A=-\frac{c}{(2n+4)\kappa} I
 \]
for some non-zero $\kappa\in\R$. Note that  $c=\frac{\Sca^{g_{o}}}{4n}$ is negative for $q\neq 0$ and positive for $q=0$.
\ee
We deduce that the quaternionic skew-Hermitian symmetric spaces $(G/L, Q, \om)$ with $G$ semisimple, which are in addition pseudo Wolf spaces, have necessarily $Q$-Hermitian Ricci tensor.
Based on \cite{AM} we may pose a kind of converse.
\bt\label{themCurvqK}
Assume that $n\geq 2$. A $4n$-dimensional quaternionic skew-Hermitian manifold $(M, Q, \om)$ with non-degenerate $Q$-Hermitian Ricci tensor  is a  quaternionic K\"ahler locally symmetric space. In particular, if $M$ is simply connected and complete, then $(M, Q, \om)$ is one of the spaces 
\[
\SU(2+p,q)/(\SU(2)\SU(p,q)\U(1))\,,\quad\ \Sl(n+1,\mathbb{H})/(\Gl(1,\mathbb{H})\Sl(n,\mathbb{H}))\,.
\]
\et
\pr
Consider a quaternionic skew-Hermitian manifold $(M, Q, \om)$ with the unique torsion-free $\SO^*(2n)\Sp(1)$-connection $\nabla^{Q, \om}$ of Proposition \ref{torsion}. To simplify the notation  let us set  $\nabla:=\nabla^{Q, \om}$ and  $\Ric:=\Ric^{\nabla}=\Ric^{Q, \om}$ for the associated Ricci tensor. 
By assumption $\Ric$ is Hermitian and non-degenerate, hence by Theorem  \ref{RicciTHM} we  have 
\begin{equation}\label{RicAsostar}
\Ric(X, Y)=(2n+4)\kappa\thinspace\om(AX, Y)\,,\quad X, Y\in\Gam(TM)\,,
\end{equation}
 for some non-zero real number $\kappa$ and some non-degenerate element $A\in\fr{so}^*(2n)$. In particular, $\Ric$ is symmetric. 
Thus, by \cite[Thm.~5.5]{AM} our quaternionic skew-Hermitian connection $\nabla$ must coincide with the Levi-Civita connection induced by the pseudo quaternionic K\"ahler metric on $M$ which is given  by $\Ric$. Next we will show that  $M$ is locally symmetric, i.e., $\nabla R=0$, where    $R=R^{Q, \om}$  is the curvature tensor of $\nabla=\nabla^{Q, \om}$, given by  (\ref{RQom}) for $A\in\fr{so}^*(2n)$ (according to Proposition \ref{adjointbundle}).
To prove this claim we use that $\nabla \Ric =0$, which holds, since $\nabla$ is the Levi-Civita connection of the metric  $\Ric$. 
By (\ref{RicAsostar}) we see that the condition $\nabla\Ric=0$ is equivalent to $\nabla A =0$,  since $\nabla\om=0$.
Let us also recall that the tensor field $\sum_a J_a\otimes J_a$ is parallel, because the connection preserves the 
quaternionic structure. In order to exploit all these  observations 
we replace in  (\ref{RQom}) the tensors $g_{a}$ by $\om (\cdot , J_a\cdot )=-\om\circ J_a$, which gives
{\small\eqnaa
R(X, Y)Z&=&\kappa\thinspace\om(X, Y)AZ+\frac{\kappa}{2}\Big(\om(X, Z)AY-\sum_{a}\om(X, J_{a}Z)J_{a}AY+\om(AY, Z)X-\sum_{a}\om(AY, J_{a}Z)J_aX\Big)\nonumber\\
&&-\frac{\kappa}{2}\Big(\om(Y, Z)AX-\sum_{a}\om(Y, J_aZ)J_{a}AX+\om(AX, Z)Y-\sum_{a}\om(AX, J_{a}Z)J_aY\Big)\nonumber\\
&&-\frac{\kappa}{2}\sum_{a}\big(\om(X, J_{a}AY)-\om(Y, J_{a}AX)\big)J_{a}Z\,.\nonumber  
\deqnaa} 

\noindent 
Then, it is clearly visible that $R$ is formed of parallel objects and hence $\nabla R=0$. If    $M$ is in addition simply connected and complete, then it must be a quaternionic K\"ahler symmetric space $M\cong G/L$. In this case,  the Lie group $G$ is necessarily semisimple (otherwise $M$ will be hyperK\"ahler, see \cite{AC}, which is Ricci flat and hence the Ricci tensor will be  degenerate).  Hence our conclusion follows by the classification of quaternionic skew-Hermitian symmetric spaces  presented in \cite{CGWPartI} and the explanation given in Example \ref{qsHpWolf}.
\pro

 \section{Bundle constructions and $\SO^*(2(n+1))$-structures}
 
 \subsection{The  Swann bundle over a quaternionic skew-Hermitian manifold}\label{knowC}
 
For the following of this section we  consider a quaternionic skew-Hermitian manifold $(M, Q, \om)$ with $\dim M=4n>4$,   together with its unique torsion-free $\SO^*(2n)\Sp(1)$-connection $\nabla:=\nabla^{Q, \om}=\nabla^{Q, \vol}$ from Theorem~\ref{connections} and Proposition \ref{torsion}.   Recall that  $\om$ is a symplectic form  which is  $Q$-Hermitian, or equivalently $\Sp(1)$-invariant (or $\fr{sp}(1)$-invariant).  Let us  denote by $\pi  : \mc{Q}\to M$  its  reduced frame bundle and by  $\pi_{S} : S\to M$  the principal $\SO(3)$-bundle of admissible frames  of the quaternionic structure $Q$.  For simplicity we will assume that the principal actions are right-actions, and in particular  we will denote the second one  by $\la_{g}(s)=sg$ for any $s=(I_1, I_2, I_3)\in S$ and $g\in\SO(3)$   (see \cite{ACDM, CH} for a more general setting). 
Note  that we have a   diffeomorphism $S\cong \mc{Q}/\SO^*(2n)$.

Next we will use the natural isomorphism $\Hn^{\times}\cong\R_{+}\times\Sp(1)$,  where $\R_{+}:=\{t\in\R : t>0\}$. The Lie algebra of  the Lie group $\Hn^{\times}$  is generated by the standard quaternionic basis $\{1, i, j, k\}$. 
In particular,  we will often denote by  $\{e_1, e_2, e_3\}$ the  basis of $\fr{so}(3)$  which corresponds to the standard basis of   $\fr{sp}(1)\cong{\rm Im}(\Hn)\cong\R^3$, under the canonical isomorphism $\fr{sp}(1)\cong\ad(\fr{sp}(1))=\fr{so}(3)$.  We also set $e_0=1\in\R\cong T_1\R_{+}$, so that we can identify $\{1, i, j, k\}$  and $\{e_0, e_1, e_2, e_3\}$. Note that  $[e_a, e_b]=2e_c$ for any cyclic permutation $(a, b, c)$ of $(1, 2, 3)$.

  The quaternionic skew-Hermitian  connection $\nabla=\nabla^{Q, \om}$   induces a principal connection  on $\mc{Q}$ and  a principal connection on $S$, with connection 1-forms
\[
 \gamma :  T\mc{Q}\to\fr{so}^*(2n)\oplus\fr{sp}(1)\,,\quad\quad   \theta :  TS\to\fr{sp}(1)\,,
\]
respectively. We have $\theta=\sum_{a=1}^{3}\theta_{a}e_{a}$, where $\theta_{a}$ are the components of $\theta$, and 
moreover we can express $\gamma$ as $\gamma=\gamma_{+}+\gamma_{-}$, with $\gamma_{+}$ taking values in $\fr{so}^*(2n)$ and $\gamma_{-}=\theta$ taking values in $\fr{sp}(1)$. 
Note that after choosing a  local frame $s=(I_1, I_2, I_3)\in S$ for $Q$ with $I_1I_2=I_3$,  the pull-backs $\pi_{S}^*\theta_{a}$ coincide with the  (local) 1-forms $\phi_{a}$ on $M$ appearing  in the relation 
\begin{equation}\label{ourtorsionfree}
 \nabla I_{a}=2(\phi_{c}\otimes I_{b}-\phi_{b}\otimes I_{c})\,,
 \end{equation}
 for  any cyclic permutation $(a, b, c)$ of $(1, 2, 3)$.
Next we write $TS=\mc{V}\oplus\mc{H}$ with $\mc{V}:=\ke\dd\pi_{S}$ and $\mc{H}:=\ke\theta$, for the the direct sum decomposition of the tangent bundle $TS$ of $S$, induced by $\theta$.

Let us also consider the (trivial) bundle 
 $S_{0}=(\wed^{4n}(T^*M)\backslash\{0\})/\Z_2$ over $M$, with fibers consisting  of volume forms at each point of  $M$ considered up to sign (it can be identified with the bundle of positive densities).  This is a  principal $\R_{+}$-bundle and we will denote by   $\pi_{S_{0}} : S_{0}\to M$ the  bundle projection, and by $\la_{0}$  the associated right principal action. 
  \bl\label{trivialS0}
 The bundle $S_{0}$ over $(M, Q, \om)$ is canonically trivial, $S_{0}=M\times\R_{+}$. in such a way that the connection 1-form
 $\theta_{0} : TS_{0}\to\R$ corresponding to the principal connection on $S_{0}$ induced by $\nabla$,  is given by $\theta_{0}=t^{-1}\dd t$, where $t$ denotes the fiber coordinate  in the trivialization. 
 \el
 \pr
 The volume form $\vol=\om^{2n}$ is a global section of $S_{0}$ and  hence induces a global trivialization of $S_{0}$.
Since the connection preserves the volume form, the second claim follows.
 \pro
Now, the product   $\rho=\la_{0}\times\la$ is also a principal action and  thus $S_{0}\times S$ is a principal bundle over $M\times M$. Its structure group is the conformal special orthogonal group in three dimensions, that is, the Lie group $\R_{+}\times\SO(3)\cong\Hn^{\times}/\Z_2$. We   consider the fiber product bundle
\[
\hat{M}=S_{0}\times_{M} S=\{(m, (v, s))\in M\times(S_{0}\times S) : m=\pi_{S_{0}}(v)=\pi_{S}(s)\}\,.
\]
In terms of the diagonal map $\vartriangle \ : M\to M\times M$ we see that $\hat{M}=\vartriangle^*\negthinspace(S_{0}\times S)$.
 Thus $\hat{M}$  is a principal $\Hn^{\times}/\Z_2$-bundle over $M$,   known as the  \textsf{Swann bundle} (see \cite[p.~222]{CH} for  details on the terminology), and we denote by $\hat{\pi} : \hat{M}\to M$ the bundle projection. 
 
 Next our goal  is to examine the geometry on the total space of  $\hat{\pi} : \hat{M}\to M$, induced  by the quaternionic skew-Hermitian structure $(Q, \om)$ downstairs on $M$. 
Let us consider the 1-form  
\begin{equation}\label{bartheta}
\hat{\theta}:=\vartriangle_{\sharp}^*\negthinspace(\theta_{0}\circ pr_{TS_{0}}, \theta\circ pr_{TS})\,,
\end{equation}
which  is a connection 1-form   on $\hat{M}$ with values in $\Hn\cong\R\oplus\fr{sp}(1)$, see also \cite{CH}.  Here  $\vartriangle_{\sharp} \ : \hat{M}\to S_{0}\times S$ is the canonical bundle morphism induced by $\vartriangle$, and $pr_{TS_{0}}$ (resp. $pr_{TS}$) is the projection from $T(S_{0}\times S)\cong TS_{0}\times TS$ onto $TS_{0}$ (resp. $TS$).
We begin by emphasizing  the following   diffeomorphisms and properties of $\hat\theta$:
\bl
Given a   quaternionic-skew Hermitian manifold $(M^{4n}, Q, \om)$, the   Swann bundle   $\hat\pi : \hat M\to M$  over $M$
satisfies 
\[
\hat{M}\cong S\times_{\SO(3)}(\R_{+}\times\SO(3))\cong \R_{+}\times S\cong \mc{Q}\times_{\SO^*(2n)\Sp(1)}(\Hn^{\times}/\Z_2)\,.
\]
 In particular, under the point of view of the diffeomorphism $\hat{M}\cong\R_{+}\times S$ the connection 1-forms $\hat\theta$ and  $\theta$ are such that
\[
\hat\theta=t^{-1}\dd t+\theta=\theta_{0}+\sum_{a=1}^{3}\theta_{a}e_{a}\,. \quad (\ast)
\]
 \el
 \pr
 The first diffeomorphism is well-known and it occurs due to the actions mentioned  above, see \cite{Swann}.  The second  one means that $\hat{M}$ is a trivial bundle over $S$.  
 The last mentioned diffeomorphism is based on the identification $S\cong \mc{Q}/\SO^*(2n)$.
Finally, the given formula $(\ast)$ is a consequence of Lemma \ref{trivialS0}
 and the triviality of $\hat{M}$ over $S$, i.e., $\hat{M}\cong\R_{+}\times S$. 
 \pro
  \be
(1) Consider the  manifold  $M=\SU(2+p,q)/(\SU(2)\SU(p,q)\U(1))$ as in Example \ref{qsHpWolf}.
The $\SO(3)$-bundle $\pi_{S} : S\to M$  has total space the $(4n+3)$-dimensional homogeneous space $\SU(2+p, q)/\SU(p, q)\U(1)$. \\
(2) Consider the    manifold $M=\SO^*(2(n+1))/\SO^*(2n)\U(1)$ as in Example \ref{notWolf}.  The $\SO(3)$-bundle $\pi_{S} : S\to M$ is given by the associated bundle 
$S=\SO^*(2(n+1))\times_{\SO^*(2n)\U(1)} \Sp(1)$, where  the action of the $\U(1)$-factor on $\Sp(1)$ is by  
left-multiplication.   
The manifolds $M$ and $S$ mentioned in this example will 
play a role in Proposition \ref{symSpaHol}. 
\ee

 \subsection{The canonical hypercomplex structure} Let $(M, Q, \om)$ be a quaternionic skew-Hermitian manifold, as above, endowed with the adapted  torsion-free connection $\nabla=\nabla^{Q, \om}$,  and let $\hat\pi : \hat{M}\to M$ be the principal $\R_{+}\times\SO(3)$-bundle   introduced in the previous subsection, with connection 1-form $\hat\theta$  given in (\ref{bartheta}). 
Relatively to $\hat\theta$ and $\hat{\pi}$ we have the direct sum decomposition 
\[
T\hat{M}=\hat{\mc{V}}\oplus\hat{\mc{H}}\,, \quad\hat{\mc{V}}:=\dd\hat\pi\,,\quad \hat{\mc{H}}:=\ke\hat\theta\,.
\]
Next we will denote by $(X)^{\hat{h}}$ the horizontal lift  of a vector field $X$ on $M$ with respect to $\hat{\theta}$.  
 Let also $\wi{U}$ be  the fundamental vector field  induced by an element  $U\in\Hn=\R\oplus\fr{so}(3)$ in the Lie algebra of the structure group of $\hat{\pi}$, that is, $\wi{U}_{u}=\frac{\dd}{\dd t}\big|_{t=0}\rho_{\exp tU}(u)\in T_{u}\hat{M}$  for any $u=(m, (v, s))\in\hat{M}$.   As in \cite{CH}, we set   
\[
Z_{a}:=\wi{e}_{a}\,,\quad a=1, 2, 3\,,\quad\quad  Z_{0}:=\wi{e}_0\,,
\] 
and consider  the triple $\hat{H}:=(\hat{I}_{1}, \hat{I}_{2}, \hat{I}_{3})$ consisting of endomorphisms $\hat{I}_{a} : T\hat{M}\to T\hat{M}$ of $T\hat{M}$ defined as follows:
\begin{equation}\label{Hthetac}
\left\{
\begin{tabular}{l}
$\hat{I}_{a}Z_{0}=-Z_a\,,\quad \hat{I}_{a}Z_{a}=Z_{0}\,,\quad\hat{I}_{a}Z_{b}=Z_{c}\,,\quad \hat{I}_{a}Z_{c}=-Z_{b}$,\\
$(\hat{I}_{a})_{u}(Y)=(I_{a}(\hat{\pi}_{*}Y))^{\bar{h}}_{u}\,,\quad \forall \ Y\in\hat{\mc{H}}_{u}$,
\end{tabular}\right.
\end{equation}
 where $u=(m, (v, s))\in\hat{M}$ and $s=(I_1, I_2, I_3)\in S$. The elements of $\hat{H}$ are almost complex structures on $\hat{M}$ satisfying  the basic quaternionic identity, $\hat{I}_{1}^2=\hat{I}_{2}^2=\hat{I}_{3}^2=\hat{I}_1\hat{I}_2\hat{I}_3=-1$,
 and thus $\hat{H}$ is an almost hypercomplex structure on $\hat{M}$, see \cite{CH} for details.

 \br
   When the initial data is (only) a quaternionic manifold $(M, Q, \nabla )$ endowed with a quaternionic 
   connection, then Lemma \ref{trivialS0} is not valid. In particular, it may not be possible to trivialize $S_0$ by a parallel section, not even locally.  
   In this setting,  a 1-parameter family $\hat{H}^{\hat{\theta}, \ccc}$ of almost hypercomplex structures on the Swann
   bundle $\hat{M}$ over $M$ was defined in  \cite[p.~216]{CH} by considering the vector field $Z_{0}^{\ccc}=\ccc\cdot\widetilde{e}_0$ for some non-zero real number $\ccc$,
   instead of our $Z_{0}$ given above (see \cite{PPS} for another approach considering a 1-parameter family of Swann bundles).  
   The integrability of the family $\hat{H}^{\hat{\theta}, \ccc}$  depends in general both  on the chosen  quaternionic connection $\nabla$ (inducing the principal connection $\hat \theta$) and on the parameter $\ccc$, see  \cite[Thm.~3.6]{CH}. However,  in our setting of quaternionic 
   skew-Hermitian structures the 
   connection $\nabla$ is unique, the almost hypercomplex structures $\hat{H}^{\hat{\theta}, \ccc}=\hat{H}^{\ccc}$ on $\hat{M}$ are all 
   hypercomplex (as explained in the proof of Theorem \ref{CHtheorem1}) and, in fact, pairwise diffeomorphic by automorphisms  
   of the trivial principal bundle $S_0$.  In particular, for $\ccc=1$ we obtain the following statement.
   \er
 
 \bt\label{CHtheorem1}
 Let $(M^{4n}, Q, \om)$ $(n>1)$ be a quaternionic skew-Hermitian manifold and let $\hat\pi :\hat{M}\to M$ be the Swann bundle over $M$.
  Then, the  almost hypercomplex structure $\hat{H}=\hat{H}^{\ccc=1}$ on $\hat{M}$ defined in {\rm(\ref{Hthetac})} is 1-integrable,  and hence a hypercomplex  structure.
  \et
  \pr
 By \cite[Thm.~3.6]{CH} we know that for a quaternionic manifold $(M, Q)$,  the almost hypercomplex structure $\hat{H}^{\hat{\theta}, \ccc}$  is independent of the choice of a quaternionic connection $\nabla$, if and only if $\ccc=-4(n+1)$ and for such  $\ccc$ the triple $\hat{H}^{\hat{\theta}, \ccc}$ is always 1-integrable.  When $\ccc\neq-4(n+1)$,    then the same theorem states that  triple $\hat{H}^{\hat{\theta}, \ccc}$ is 1-integrable if and only if $(\Ric^{\nabla})^{a}$ is $Q$-Hermitian, where $(\Ric^{\nabla})^{a}$ denotes the skew-symmetric part of the Ricci tensor $\Ric^{\nabla}$.   
 
 In our case,   the base space $(M, Q, \om)$ of $\hat{\pi} : \hat{M}\to M$ admits a unique (torsion-free) adapted connection $\nabla$ (quaternionic-skew Hermitian connection), given by $\nabla^{Q, \om}$ (recall by   Section \ref{AdaConnec} that the first prolongation of such geometric structures vanishes). 
Therefore,  we do {\it not} have the choice of a quaternionic connection mentioned for the quaternionic case, but a unique one: $\nabla=\nabla^{Q, \om}$. On the other hand, by Corollary \ref{QHermitianRic} the 
 Ricci tensor of $\nabla^{Q, \om}$ is always symmetric.
 Thus $\Ric^{\nabla}$ has zero skew-symmetric part and trivially it is $Q$-Hermitian. The claim now follows by \cite[Thm.~3.6]{CH}.
  \pro

 \subsection{$\SO^*(2(n+1))$-structures on $\hat{M}$}
In the sequel  we focus on the description of $\SO^*(2(n+1))$-structures on the total space $\hat{M}$ of the Swann bundle $\hat{\pi} : \hat{M}\to M$ over a quaternionic skew-Hermitian manifold $(M^{4n}, Q, \om)$ $(n>1)$. We will show that such structures are  in general non-integrable and  of type $\mc{X}_{3457}$, in terms of the  intrinsic torsion characterization of almost hypercomplex skew-Hermitian structures given in  \cite{CGWPartI}.  In other words, below we will prove that  except of the  flat case,  $\hat{M}$ comes only with pairs $(\hat{H}, \upomega)$, where $\upomega$ is  an $\hat{H}$-Hermitian almost symplectic form, but not symplectic (i.e., not closed). We begin with a description of   $\hat{H}$-Hermitian 2-forms on the  horizontal part of the decomposition $T\hat{M}=\hat{\mc{H}}\oplus\hat{\mc{V}}$. These are induced by $\om$, as the following proposition verifies.

\bp\label{hatomega}
 Let $(M, Q, \om)$ and $\hat{\pi} : \hat{M}\to M$ be as above, and set $G:=\R_{+}\times\SO(3)$.\\
$(1)$ The pullback 
\[
\hat{\om}:=\hat{\pi}^*(\om)\in\Om^{2}(\hat{M})
\]
 defines   a horizontal,  $G$-invariant,    closed 2-form   on $\hat{M}$, satisfying $\mc{L}_{Z}\hat\om=0$ for any vertical vector $Z\in\hat{\mc{V}}$.\\
 $(2)$  The restriction $\hat\om_{\hat{\mc{H}}}:=\hat\om\big|_{\hat{\mc{H}}\times\hat{\mc{H}}}$ of $\hat{\om}$ to the horizontal subspace $\hat{\mc{H}}=\ke\hat\theta\subset T\hat{M}$ is non-degenerate, closed and $\hat{H}$-Hermitian. 
Note also that since $\hat\om$ is horizontal, the differential of $\hat\om_{\hat{\mc{H}}}$ is defined as a section of $\bigwedge^3  {\hat{\mc{H}}}$ and coincides with the restriction of the differential of $\hat\om$.
\ep
\pr
(1) The closedness is obvious. Consider the decomposition $T\hat{M}=\hat{\mc{V}}\oplus\hat{\mc{H}}$ given above. The fundamental vector fields $\{Z_{0}, Z_{a} : a=1, 2, 3\}$ span $\hat{\mc{V}}$ and by the definition of $\hat\om$ it is direct   that $\hat{\om}(Z_{0}, Z_{a})=\hat{\om}(Z_{a}, Z_{b})=0$.  
In addition,  for any $Y\in\hat{\mc{H}}$ and $Z\in\hat{\mc{V}}$ we have  $\hat{\om}(Y, Z)=\om(\hat\pi_{*}Y, \hat\pi_*Z)=0$, since $Z\in\ke\hat\pi_*$.  Therefore, for any $Z\in\hat{\mc{V}}$ we obtain the relation  $\iota_{Z}\hat\om=\hat\om(Z, -)=0$ and $\hat{\om}$ annihilates the vertical vectors, so it is horizontal. The relation  $\mc{L}_{Z}\hat\om=0$ follows now by the Cartan's formula. \\
 Let us now consider the action of the structural group $G$ on $\hat{M}$, denoted  by $\rho$. We have
  \[
\rho_{(t, g)}\big(m, (v, s)\big):=\big(m, (\la_0(t, v), \la(g, s))\big)=\big(m, (v t, s g)\big)\,, \quad (m, (v, s))\in\hat{M} 
  \]
  for any $t\in\R_{+}$ and $g\in\SO(3)$. Since $\hat{\pi} : \hat{M}\to M$ is a principal $G$-bundle we have that $\hat\pi\circ\rho_{(t, g)}=\hat\pi$ for any $(t, g)\in G$ and hence    $\rho_{(t, g)}^*\hat\om=\rho_{(t, g)}^*(\hat\pi^*(\om))=(\hat\pi\circ\rho_{(t, g)})^*\om=\hat\pi^*\om=\hat\om$. 
  This shows that $\hat\om$ is $G$-invariant.\\
(2) To prove the statement for  non-degeneracy consider the equation  $\hat\om(Y, Z)=0$ with   $Y, Z\in\Gam(\hat{\mc{H}})$. This  means that  $\om(\hat\pi_*Y, \hat\pi_*Z)=0$. Since $\hat\pi$ is a submersion and $\om$ is non-degenerate, we deduce that $\hat\pi_*Y=0$, which implies that $Y=0$, since  by assumption $Y$ is not vertical. Finally, let $Y, Z\in\hat{\mc{H}}_{u}$ be two horizontal vectors at $u=(m, (v, s))\in\hat{M}$. Then we see that
\eqna
\hat{\om}_{u}\big((\hat{I}_{a})_{u}Y, (\hat{I}_{a})_{u}Z\big)&=&\om_{\hat{\pi}(u)}\Big(\hat{\pi}_{*}\big((\hat{I}_{a})_{u}Y\big), \hat{\pi}_{*}\big((\hat{I}_{a})_{u}Z\big)\Big)\\
&\overset{(\ref{Hthetac})}{=}&\om_{\hat{\pi}(u)}\Big(\hat{\pi}_{*}\big((I_{a}(\hat{\pi}_{*}Y))^{\bar{h}}_{u}\big), \hat{\pi}_{*}\big((I_{a}(\hat{\pi}_{*}Z))^{\bar{h}}_{u}\big)\Big)\\
&=&\om_{\hat{\pi}(u)}\big(I_{a}(\hat{\pi}_{*}Y), I_{a}(\hat{\pi}_{*}Z)\big)\\
&=&\om_{\hat{\pi}(u)}\big(\hat{\pi}_{*}Y, \hat{\pi}_{*}Z\big)=\hat{\om}_{u}(Y, Z)
\deqna
for any $a=1, 2, 3$, where the third equality follows by the definition of the horizontal lift, and the fourth  since $\om$ is $Q$-Hermitian and $s=(I_1, I_2, I_3)$ is an admissible basis of $Q$.
\pro
As a first immediate    consequence of the   above result we have that
\bc\label{SO*parallelpart}
For any $u\in \hat{M}$, the triple $(\hat{\mc{H}}_{u}:=\ke\hat{\theta}_{u}\subset T_{u}\hat{M},  \hat{H}_{\hat{\mc{H}}_{u}}, \hat{\om}_{\hat{\mc{H}}_{u}})$ is a hypercomplex skew-Hermitian vector space, where  
\[
\hat{H}_{\hat{\mc{H}}_{u}}:=\{\hat{I}_{a}\in H : (\hat{I}_{a})_{u}(Y)=(I_{a}(\hat{\pi}_{*}Y))^{\bar{h}}_{u}\,, \ \text{for all} \ Y\in\hat{\mc{H}}_{u}\}\,.
\]
\ec
About the 2-form $\hat\om$, by Proposition \ref{hatomega} we also deduce the following. 
\bc\label{presymplectic}
The 2-form $\hat\om$ is a basic differential form on $\hat{M}$ which is  presymplectic. Its kernel $\ke\hat\om=\{X\in\Gam(T\hat{M}) : \iota_{X}\hat\om=0\}$ is the linear span $\lan Z_{0}, Z_{1}, Z_{2}, Z_{3}\ran = \hat{\mc{V}}$, which is an integrable distribution.
\ec

Let us now look for  $\hat{H}$-Hermitian 2-forms on the vertical part.
For this it is useful to consider the $\Hn$-valued 1-form $\al$ on $\hat{M}$ given by
\[
\al:= \frac{1}{t}\hat{\theta}= \frac1{t^2}\dd t+ \frac{1}{t}\theta\in\Om^{1}(\hat{M}; \Hn )=\Om^{1}(\hat{M}; \mathbb{R} \oplus \mathfrak{sp}(1))\,.
\]
Obviously, we have $\al=\alpha_0+\al_{1}e_1+\al_{2}e_2+\al_{3}e_3$ for some real-valued 1-forms $\al_a$ on $\hat{M}$. Moreover,
by expressing the connection 1-form $\theta$ as $\theta=\sum_{a}^3 \theta_ae_a$, we see that
\begin{equation}\label{ddaleijk}
\al_{0}= \frac1{t^2}\dd t\,, \quad \al_{a}=\frac{1}{t}\theta_{a}\,,\quad a=1,2,3.
\end{equation}
\br
Note that $\al$  has some similarity with the $\Hn$-valued 1-form introduced on $\hat{M}$ over a (pseudo) quaternionic-K\"ahler manifold $(M^{4n}, Q, g)$ by Swann \cite{Swann}.
\er

\bp\label{32forms}At any point $p\in \hat{M}$, 
the space of all $\hat{H}$-Hermitian 2-forms    $\beta$ with $\beta|_{\hat{\mc{H}}_p\times\hat{\mc{H}}_p}=0$, 
$\beta|_{\hat{\mc{H}}_p\times \hat{\mc{V}}_p}=0$ is 3-dimensional, generated by the  2-forms
\[ \beta_a = \alpha_0 \wedge \alpha_a + \alpha_b\wedge \alpha_c \]
for any  cyclic permutation $(a,b,c)$ of $(1,2,3)$. 
\ep 
\pr
Since the form $\beta$ vanishes on the horizontal distribution, it reduces to a 
vector-valued 2-form on the 4-dimensional vertical distribution. 
Using the fact that the four components of the connection form $\hat\theta$ are dual to the frame $(Z_i)_{i=0,1,2,3}$ of the vertical space of 
$\hat M$ and that $\alpha =\frac1t \hat\theta$ we see that the structures $\hat{I}_a$, $a =1,2,3$, act as follows on the $\alpha_i$, $i=0, \ldots, 3$: 
\[ \hat{I}_a^*\alpha_0 = \alpha_a,\qquad \hat{I}_a^*\alpha_b = - \alpha_c\,,\]
for any cyclic permutation.  From these formulas one can immediately check that the forms $\beta_a$ are invariant under
the hypercomplex structure $\hat{H}$.  The remaining statement follows from the fact that the space of 
$\mathfrak{sp}(1)$-invariant 2-forms on $V=\mathbb{R}^4$ is 3-dimensional (which corresponds to one summand in 
the decomposition $\wedge^2V^* \cong \mathfrak{so}(4) = \mathfrak{sp}(1) \oplus \mathfrak{sp}(1)$). 
\pro

 \bt\label{mainthmS}
 Let $\hat\pi : \hat M\to M$ be  the Swann bundle over a quaternionic-skew Hermitian manifold $(M^{4n}, Q, \om)$ $(n>1)$, the latter endowed with the unique torsion-free quaternionic skew-Hermitian
 connection $\nabla=\nabla^{Q, \om}$. Let $R^{Q, \om}$ be the curvature tensor induced by $\nabla^{Q, \om}$. If $R^{Q, \om}\neq 0$, then there is no non-degenerate closed vertical  
 $\hat{H}$-Hermitian 2-form $\beta$ on the Swann bundle $\hat{M}$, i.e., scalar 2-form, such that $\beta|_{\hat{\mc{H}}\times\hat{\mc{H}}}=0$, $\beta|_{\hat{\mc{H}}\times\hat{\mc{V}}}=0$, $\beta|_{\hat{\mc{V}}\times\hat{\mc{V}}}$ non-degenerate and $\dd\beta=0$. 
In particular, it is not possible to extend the canonical closed $\hat{H}$-Hermitian horizontal $2$-form $\hat{\omega}=\hat\pi^*(\om)$ to a $\hat{H}$-Hermitian symplectic form by adding a vertical form.
 \et
 \pr
 By the previous proposition, any vertical 
 $\hat{H}$-Hermitian 2-form $\beta$ is of the form $\beta=\sum f_a \beta_a$, $f_{a}\in C^{\infty}(\hat{M})$, $a=1, 2, 3$. To  compute its differential we need to understand $\dd\beta_a$ for $a=1, 2, 3$, and for this it is necessary to compute $\dd\al_i$ for all $i=0,1,2,3$. At this point we rely on the identity
 \begin{equation}\label{ddal2}
 \dd \al= \frac{-1}{t^2}\dd t\wedge\theta+\frac1t\dd\theta=-\al_0\wedge\theta+\frac1t\Omega-\frac1t\theta\wedge\theta\,,
 \end{equation}
 where $\Om=\dd\theta+\theta\wedge\theta$ is the part of the curvature with values in $\fr{sp}(1)$ (and so the curvature induced by   the connection 1-form $\theta$, which recall that is a part of $\gamma$).  Hence, having  in mind the relations $\al_a=\frac1t\theta_a$ for $a=1, 2, 3$, $\theta=\sum_{a=1}^{3}\theta_{a}e_a$ and
 $\Om=\sum_{a=1}^{3}\Om_{a}e_a$, an application of (\ref{ddal2}) gives
 \begin{eqnarray}
 \dd \al&=&-\al_{0}\wedge(\theta_{1}e_1+\theta_{2}e_2+\theta_{3}e_3)+ \frac1t(\Om_{1}e_1+\Om_{2}e_2+\Om_{3}e_3)\nonumber\\
&-&\frac{1}{2}[\al_{1}e_1+\al_{2}e_2+\al_{3}e_3,\theta_{1}e_1+\theta_{2}e_2+\theta_{3}e_3]\nonumber\\
 &=&\sum_{cycl}(-\al_0\wedge\theta_a-2\al_b\wedge\theta_c+ \frac1t\Om_a)e_a\,. \label{ddallong}
 \end{eqnarray}
On the other hand we have
 \begin{equation}\label{ddal1}
 \dd\al=\dd\al_{0}+\dd\al_{1}e_1+\dd\al_{2}e_2+\dd\al_{3}e_3
 \end{equation}
and comparing these two equations we obtain 
 \begin{eqnarray*}
 \dd\al_{a}&=&-\al_0\wedge\theta_a-2\al_b\wedge\theta_c+ \frac1t\Om_a\,,
\end{eqnarray*}
for any  cyclic permutation $(a,b,c)$ of $(1,2,3)$, in addition to the trivial equation $d\alpha_0=0$. 
Hence we can proceed with $\dd\beta$. 
A direct computation based on the definition of $\beta_{a}$ $(a=1, 2, 3)$ and of the given formulas for the differentials of $\al_i$, $(i=0,1, 2, 3)$ shows that
\begin{eqnarray*}
\dd\beta_a&=&-\al_0\wedge\dd\al_a+\dd\al_b\wedge\al_c-\al_b\wedge\dd\al_c\\
&=&-\al_0\wedge(-2\al_b\wedge\theta_c+ \frac{1}{t}\Om_a)+
(-\al_0\wedge\theta_b+\frac1t\Om_b)\wedge\al_c-\al_b\wedge(-\al_0\wedge\theta_c+\frac1t\Om_c)\\
&=&\frac1t(-\al_0\wedge\Om_a +\Om_b\wedge\al_c-\al_b\wedge\Om_c)\,.
\end{eqnarray*}
Hence it follows that
 \begin{eqnarray}
 \dd\beta&=&\sum_{a=1}^3(\dd f_a\wedge\beta_a+f_a\dd\beta_a) \nonumber\\
 &=&\sum_{a=1}^{3}\dd f_{a}\wedge\beta_{a}-\frac1t\al_{0}\wedge(f_{1}\Om_{1}+f_{2}\Om_{2}+f_{3}\Om_{3}) \label{dbeta}\\
 &&+\frac1t\sum_{cycl} f_a(\Omega_b\wedge \alpha_c-\Omega_c\wedge \alpha_b)\,. \nonumber
 \end{eqnarray}
 From this relation we deduce that $\dd\beta=0$ if and only if the following conditions hold:
\[
\sum_{a=1}^{3}\left( \dd f_{a}\wedge\beta_{a}\right) =0\,,\quad \sum_{a=1}^{3}f_{a}\Om_{a}=0\,, \]
\[ 
f_2\Om_3-f_3\Om_2=0\,,\quad 
-f_1\Om_3+f_3\Om_1=0\,,\quad
f_1\Om_2-f_2\Om_1=0\,.
\]
Assume now that $\beta$ is non-degenerate, and let $x\in\hat{M}$ be a point of the Swann bundle such that $R^{Q, \om}(\hat{\pi}(x))\neq 0$. By the first assumption there is $a$ such that $f_a(x)\neq 0$. From the last three equations we can express $\Om_b|_x$ for $b\neq a$ as functional multiple of $\Om_a|_x$. In particular, the second equation becomes $(f_1(x)^2+f_2(x)^2+f_3(x)^2)f_a^{-1}(x)\Om_a|_x=0$.  Thus $\Om_a|_x=0$, which implies $\Om_\ell|_x=0$, for $\ell=b$ and $\ell=c$ as well. This contradicts the assumption $R^{Q, \om}(\hat{\pi}(x))\neq 0$.  Thus there is no non-degenerate $\beta$ such that $\dd\beta=0$ if $R^{Q, \om}\neq 0$.  
 \pro
 
 \br
Let us point out that  it is not accidental that we  do {\it not} consider 2-forms defined in $\hat{\mc{H}}\times\hat{\mc{V}}$, but only horizontal and vertical ones. 
This is because there are no non-zero such $\SO^*(2n)$-invariant 2-forms for $n>1$. Indeed, recall that  by representation theory of $\SO^*(2n)$ and under the $\SO^*(2n)$-action we can identify  $\E=[\E\Hh]$, where $\E=\C^{2n}$ is the standard representation of $\SO^*(2n)$   (see \cite{CGWPartI}). Since  $\hat{\mc{H}}$ coincides pointwise with $[\E\Hh]$, our claim follows.  
\er

\bc \label{Cor_main}
Let $\hat\pi : \hat M\to M$ be the Swann bundle over a quaternionic-skew Hermitian manifold $(M^{4n}, Q, \om)$ $(n>1)$, the latter endowed with the torsion-free quaternionic skew-Hermitian
 connection $\nabla=\nabla^{Q, \om}$. Set $\tilde\om:=\hat\om+\beta$,  where $\hat\om=\hat\pi^*(\om)$, $\beta=\sum_{a=1}^{3}f_{a}\beta_{a}$ for some   smooth functions $f_{a}$  on $\hat{M}$, such that the vector-valued function $(f_a)$ is nowhere vanishing and $\beta_{a}$ $(a=1, 2, 3)$ are as in Proposition \ref{32forms}. 
Then, the  the following statements hold. 
\begin{enumerate}
\item The
pair $(\hat{H}, \tilde\om)$ defines an $\SO^*(2(n+1))$-structure on the total space $\hat{M}$ of the Swann bundle over $M$, which is in general of type $\mc{X}_{3457}$. 
\item 
The  $\SO^*(2(n+1))$-structure $(\hat{H},\tilde\om)$ has always nontrivial torsion in component $\mc{X}_{34}$, provided that $R^{Q, \om}\neq 0$. 
\item
Each non-zero vector $(v_a)\in \mathbb{R}^3$ defines a canonical $\SO^*(2(n+1))$-structure 
of type $\mc{X}_{3457}$ by choosing in the above construction the constant functions $f_a=v_a$, $a=1,2,3$.
\end{enumerate}
\ec
\pr
By the definition of $\tilde\omega$ we have
\begin{equation}\label{tildeomeg}
\tilde\om(X, Y):=\left\{
\begin{tabular}{l l}
$\hat\om(X, Y)$\,, & $\text{if} \  \ X, Y\in\hat{\mc{H}}$\,,\\
$0$\,, & $\text{if} \ \ X\in\hat{\mc{V}}, \ Y\in\hat{\mc{H}}$\,,\\
$\beta(X, Y)$\,, & $\text{if}  \ \ X, Y\in\hat{\mc{V}}$\,.
\end{tabular}\right. 
\end{equation}
Now,  by (2) in  Proposition \ref{hatomega} we deduce that  $\hat\om^{2n}\neq 0$, that is,  $\hat\om$ is non-degenerate on $\hat{\mc{H}}$. 
Also, since   we assume that $(f_{a})$ is nowhere vanishing, the 2-form  $\beta=\sum_{a=1}^{3}f_{a}\beta_a$ is nowhere vanishing and hence non-degenerate on $\hat{\mc{V}}$. 
All together and having in mind the relation (\ref{tildeomeg}) a direct computation shows that  $\tilde\omega^{2(n+1)}$ is a volume form on $\hat{M}$, that is, 
\[
\tilde\om^{2(n+1)}=(\hat\om+\beta)^{2(n+1)}=\beta^2 \wedge \hat\om^{2n}\neq 0
\]
(note that terms  containing $\beta$ to a power greater then $2$, or $\hat\om$ to a power greater than $2n$,  are  zero).
Thus $\tilde\om$ is non-degenerate.  
Combining now  Propositions \ref{hatomega} and \ref{32forms} we deduce that the pair $(\hat{H}, \tilde\om)$   defines an  $\SO^*(2(n+1))$-structure on $\hat{M}$.  As $\hat{H}$ is 1-integrable, we know that this structure is in general of type $\mc{X}_{3457}$, see \cite{CGWPartII}. The failure of $\dd\tilde\om=\dd(\hat\om+\beta)=\dd\beta$ to vanish, contributes with the torsion components $\mc{X}_{34}$, see \cite[Theorem 1.5]{CGWPartII}. 
This completes the proof.
\pro

We know from  \cite{CGWPartII} that $\SO^*(2(n+1))$-structure of type $\mc{X}_{457}$ are  locally conformally symplectic, i.e.,  locally there is functional multiple of the scalar 2-form which is a closed scalar 2-form. Let us investigate the possibility for the $\SO^*(2(n+1))$-structure $(\hat{H},\tilde\om)$ defined above on the Swann bundle, to be of the type $\mc{X}_{457}$. In  \cite{CGWPartII},  it was also shown  that $\dd \tilde\om$ will belong to the component  $\mc{X}_{457}$, if and only if $\dd\tilde\om=\dd\beta=\tau\wedge \tilde\om$ for a 1-form $\tau$ on $\hat M$. Note that from $\dd^2=0$ and non-degeneracy of $\tilde\om$, it follows that $\dd \tau=0.$ We get the following characterization  in this case.

\bp\label{symSpaHol}
The $\SO^*(2(n+1))$-structure $(\hat{H},\tilde\om)$ on the Swann bundle  $\hat{M}$ over a non-flat torsion-free quaternionic skew-Hermitian manifold $(M,Q,\om)$ is of the type $\mc{X}_{457}$ if and only if $M$ is locally equivalent to the symmetric space $\SO^*(2n+2)/\SO^*(2n)\U(1)$ and
\[
f_a=-\frac{4\,c_a\,n \, t^4}{cc_1(-h_0^2-h_1^2+h_2^2+h_3^2)+2cc_2(h_0h_3-h_1h_2)-2cc_3(h_0h_2+h_1h_3)+2cc_4t^4}
\]
for constants $c_1,c_2, c_3, c_4$ with $(c_1,c_2,c_3)\in \mathbb{R}^3 \setminus \{ 0\}$ and (local) coordinates $h_0+h_1i+h_2j+h_3k$ on the fiber $\Hn^{\times}/\mathbb{Z}_2$,
where $t=\sqrt{h_0^2+h_1^2+h_2^2+h_3^2}$ and    $c$ is the Einstein constant of $M$ described in  Example \ref{notWolf}.
\ep

\pr
Suppose $\dd\tilde\om=\dd\beta=\tau\wedge \tilde\om$ for a closed 1-form $\tau$ on $\hat{M}.$ 
Then, from the equation \eqref{dbeta} we obtain the conditions
\begin{eqnarray*}
\tau\wedge \beta&=&\sum_{a=1}^{3}\dd f_{a}\wedge\beta_{a}\,,\\
\tau\wedge \hat\om&=&- \frac1t\al_{0}\wedge(f_{1}\Om_{1}+f_{2}\Om_{2}+f_{3}\Om_{3}) + {\frac1t}\sum_{cycl} f_a(\Omega_b\wedge \alpha_c-\Omega_c\wedge \alpha_b)\,.
\end{eqnarray*}
By the non-degeneracy  of $\beta_a$ on $\hat{\mc{V}}$ we see that $f_a\tau=\dd f_{a}$ for $a=1,2,3,$ i.e., there is a function $f$ on $\hat{M}$ and constants $c_a\in \R$ such that $f_a=c_a\exp(f)$ and $\tau=\dd f.$ As  $\alpha_0,\dots, \alpha_3$ are linearly independent, it follows that 
\begin{eqnarray*}
\exp(f)(c_1\Om_{1}+c_{2}\Om_{2}+c_{3}\Om_{3})&=&s_0 \hat\om\,,\\
 \exp(f)(c_2\Om_3-c_3\Om_2)&=&s_1 \hat\om\,,\\
  \exp(f)(-c_1\Om_3+c_3\Om_1)&=&s_2 \hat\om,\\ 
  \exp(f)(c_1\Om_2-c_2\Om_1)&=&s_3 \hat\om\,,\\
  \tau&=&-\frac1ts_0\al_{0}+\frac1t\sum_{a=1}^{3}s_a\al_{a}\,,
  \end{eqnarray*}
  where   $s_0,\dots,s_3$ are some functions on $\hat{M}.$ This yields  that 
  \[
  \Om_{a}=r_a\hat\om\,,\quad a=1, 2, 3
  \]
  for some functions $r_1, r_2, r_3$ on $\hat{M}$. From Lemma \ref{hequivlem} and Corollary \ref{curvatureRQom}, we know which parts of $R^{Q, \om}=R_A$ contribute to $\Om_{a}$. 
 More precisely, along the section $p\mapsto 1|_{p}$  of $\hat{M}=\R_+\times S \to M$ the elements $e_1, e_2, e_3$
 (viewed as endomorphisms), pull back to $J_1, J_2, J_3$.  Moreover, we find that 
\[ \sum \Omega_a(x,y)J_a = 2\kappa \omega(x,y) A^{\mathfrak{sp}(1)} + \kappa \sum g_a(x, A^{\mathfrak{so}^*(2n)})J_a, \quad x,y\in T_pM,\]
and we can conclude that $\Om_{a}$ is proportional to $\hat\om$ if and only if $A$ is section with values in $\sp(1).$ Notice here we are 
identifying basic differential forms on $\hat M$ with forms on $M$ and, in particular, $\hat \omega$ with $\omega$.  
Now, along the previous section of $\hat M$ we have  $A=\frac{r_1e_1+r_2e_2+r_3e_3}{2\kappa}$ and the functions $r_1,r_2,r_3$ are completely determined by $A$ away from this section, using the right action of $\R_+\times \SO(3)$. Then $s_0=\exp(f)(c_1r_1+c_2r_2+c_3r_3),s_1=\exp(f)(c_2r_3-c_3r_2),s_2=\exp(f)(-c_1r_3+c_3r_1),s_3=\exp(f)(c_1r_2-c_2r_1)$ and 
\begin{eqnarray}\label{x4eqat}
\tau&=&\dd f=\frac1t\exp(f)(-(c_1r_1+c_2r_2+c_3r_3)\al_{0}+(c_2r_3-c_3r_2)\al_{1}\\
&&+(-c_1r_3+c_3r_1)\al_{2}+(c_1r_2-c_2r_1)\al_{3})\,.\nonumber
\end{eqnarray}
Now, $0=\dd \tau$ can be decomposed  according to $T\hat M = \hat{\mc{H}}\oplus \hat{\mc{V}}$. The component in 
$\hat{\mc{V}}^0 \wedge \hat{\mc{H}}^0$, where the label $0$ stands for the annihilator in $T^*\hat M$, is determined by 
restriction of 
\begin{eqnarray*}
&&\exp(f)(-(c_1\dd r_1+c_2\dd r_2+c_3\dd r_3)\wedge \al_{0}+(c_2\dd r_3-c_3\dd r_2)\wedge \al_{1}\\
&&+(-c_1\dd r_3+c_3\dd r_1)\wedge \al_{2}+(c_1\dd r_2-c_2\dd r_1)\wedge \al_{3})\,.
\end{eqnarray*}
Using the same arguments as in the proof of Theorem \ref{mainthmS}, we deduce that   $r_1,r_2,r_3$ are constant on the section $M\to \R_+\times S=\hat{M}$. 
Thus, $(M,\nabla^{Q,\omega})$ is  a locally symmetric space.  Moreover, from the formula given in Corollary \ref{curvatureRQom} 
and the Ambrose-Singer theorem we deduce that it has holonomy contained in  $\SO^*(2n)\U(1)$.  
From the classification of (non-flat) symmetric spaces it is known that there exists only one series of symmetric spaces  with such holonomy, namely the family  $M=\SO^*(2n+2)/\SO^*(2n)\U(1)$.

Let us now solve the equation \eqref{x4eqat} for this family, 
where we identify the generator of $\fr{u}(1)$ with $e_1$ on the previous section $M\to \R_+\times S=\hat{M}$.
Recall  from Example \ref{notWolf} that for this symmetric space  the constant element $A\in\fr{sp}(1)$ along
this section has the form $  A=-\frac{c}{4n\kappa}e_1$, 
where $c$ is the Einstein constant.  By introducing the  local coordinates $h_0+h_1i+h_2j+h_3k$ 
on the fiber $\Hn^{\times}/\mathbb{Z}_2$, where the image of the section above corresponds to $(1, 0, 0, 0)$,  we have
\begin{eqnarray}
\nonumber t&=&\sqrt{h_0^2+h_1^2+h_2^2+h_3^2}\,,\\
\nonumber\dd t&=&(\sqrt{h_0^2+h_1^2+h_2^2+h_3^2})^{-1}(h_0\dd h_0+h_1 \dd h_1+h_2\dd h_2+h_3 \dd h_3)\,,\\ 
\nonumber\theta_1|_{\mathrm{fiber}}&=&(h_0^2+h_1^2+h_2^2+h_3^2)^{-1}(-h_1 \dd h_0+h_0 \dd h_1+h_3 \dd h_2-h_2 \dd h_3)\,,\\
\nonumber\theta_2|_{\mathrm{fiber}}&=&(h_0^2+h_1^2+h_2^2+h_3^2)^{-1}(-h_2 \dd h_0-h_3 \dd h_1+h_0 \dd h_2+h_1 \dd h_3)\,,\\ 
\label{array:eq} \theta_3|_{\mathrm{fiber}}&=&(h_0^2+h_1^2+h_2^2+h_3^2)^{-1}(-h_3 \dd h_0+h_2 \dd h_1-h_1 \dd h_2+h_0 \dd h_3)\,,
\end{eqnarray}
where the restriction is to the fiber of $\hat{M} \to M$ and we have used that, on the fiber, the 1-form $\hat\theta = t^{-1}dt + \theta$ coincides with 
the left-invariant Maurer-Cartan form.
Using now the principal (right) action we can determine  $r=r(h)=r_1(h)i +r_2(h)j+ r(h)_3k$ explicitly by computing  the expression
\[ r(h) = -\frac{c}{2n}\mathrm{Ad}_{h}^{-1}(e_1)=-\frac{c}{2n}h^{-1}ih\,.\]
The result is:
\[
r_1=-\frac{c}{2n}\cdot\frac{h_0^2+h_1^2-h_2^2-h_3^2}{t^2}\,, \quad  r_2=\frac{c}{2n}\cdot\frac{2(h_0h_3-h_1h_2)}{t^2}\,, \quad  r_3=-\frac{c}{2n}\cdot\frac{2(h_0h_2+h_1h_3)}{t^2}\,.
\]
 The equation  \eqref{x4eqat}  reduces now to the standard problem of finding a primitive for a closed one-form.
Using the command ${\tt{pdsolve}}$ in Maple\footnote{See the maple worksheet at 
\tt{http://www.math.uni-hamburg.de/home/cortes/CCG.zip}} 
we find the general solution $f$, determined by 
\[  
\exp (f)=-\frac{4n \, t^4}{cc_1(-h_0^2-h_1^2+h_2^2+h_3^2)+2cc_2(h_0h_3-h_1h_2)-2cc_3(h_0h_2+h_1h_3)+2cc_4t^4}\,.
\] 
This implies the claimed formula for the functions $f_a=c_a\exp (f)$.
\pro

\br
Observe that the functions $f_a$ are  well-defined as rational functions on the fiber.
\er

 \subsection{$\SO^*(2(n+1))$-structures on $\hat{M}$ over the  flat  model}

We will now consider a quaternionic-skew Hermitian manifold $(M^{4n}, Q, \om)$ $(n>1)$ 
endowed with the torsion-free quaternionic skew-Hermitian connection $\nabla=\nabla^{Q, \om}$ such that $R^{Q, \om}=0.$    
In this  case the  exponential coordinates   induced by $\nabla$ identify open subsets of $M$ with open subsets in the flat model $([\E\Hh],Q_0,\omega_0).$ 
Moreover, the transition function between two such exponential coordinates is a
composition of  translations and the linear action of $\SO^*(2n)$ on $[\E\Hh]$. 
Therefore,  in this case the Swann bundle $\hat{M}$ over $M$ is trivial, i.e., $\hat{M}=M\times \Hn^{\times}/\Z_2$ 
and moreover we can identify  the horizontal bundle with  the pull-back of $TM$, i.e., $\hat{\mc{H}}= \pi^*TM$.  
This allows us to reformulate the condition $\dd \beta=0$ as follows:

\bl\label{lemmaFlat}
Under the above assumption,  we continue using the map  $\R^4\to \Hn^{\times}/\Z_2$, 
\[
(h_0,h_1,h_2,h_3)\mapsto h_0+h_1i+h_2j+h_3k\pmod{\Z_2}\,.
\]
Then
\[
\beta=   F_1(\dd h_0\wedge \dd h_1+\dd h_2\wedge \dd h_3)+F_2(\dd h_0\wedge \dd h_2-\dd h_1\wedge \dd h_3)+F_3(\dd h_0\wedge \dd h_3+\dd h_1\wedge \dd h_2)\,,
\]
for some functions $F_a$ on $M$. 
Thus, $\dd \beta=0$ if and only if $F_a$ are constant on $M$ and satisfy the following PDE's:
\begin{eqnarray*}
F_{1,0}-F_{3,2}+F_{2,3}&=&0,\\
F_{1,1}+F_{2,2}+F_{3,3}&=&0,\\
F_{1,3}-F_{2,0}-F_{3,1}&=&0,\\
F_{1,2}-F_{2,1}+F_{3,0}&=&0,
\end{eqnarray*}
where $F_{a,b}$ denotes the partial derivative of $F_a$ with respect to $h_b.$
\el
\pr
As before,  we can use the equations given in \eqref{array:eq}. Now, to obtain the claimed formula for $\beta$ one needs to consider the following combinations of the functions $F_a$ substituted for the $f_a$ functions
\begin{eqnarray*}
f_1&=&(h_0^2+h_1^2+h_2^2+h_3^2)\left((h_0^2+h_1^2-h_2^2-h_3^2)F_1+2(h_0h_3+h_1h_2)F_2-2(h_0h_2-h_1h_3)F_3\right),\\
f_2&=&(h_0^2+h_1^2+h_2^2+h_3^2)\left(-2(h_0h_3-h_1h_2)F_1+(h_0^2-h_1^2+h_2^2-h_3^2)F_2+2(h_0h_1+h_2h_3)F_3\right),\\
f_3&=&(h_0^2+h_1^2+h_2^2+h_3^2)\left(2(h_0h_2+h_1h_3)F_1-2(h_0h_1-h_2h_3)F_2+(h_0^2-h_1^2-h_2^2+h_3^2)F_3\right).
\end{eqnarray*}

This immediately leads to the claimed formula for $\beta$ and moreover we see that
 \[
\dd \beta=\dd F_1\wedge (\dd h_0\wedge \dd h_1+\dd h_2\wedge \dd h_3)+\dd F_2\wedge (\dd h_0\wedge \dd h_2-\dd h_1\wedge \dd h_3)+\dd F_3 \wedge(\dd h_0\wedge \dd h_3+\dd h_1\wedge \dd h_2)\,.
\]
If $\dd \beta=0$, then clearly $F_a$ has to be constant on $M$, while expressing $\dd F_a=\sum_{b=0}^3F_{a,b}\dd h_b$ leads to the claimed PDE's.
\pro

The PDE system from the above lemma can be solved using Maple, however the general solution involves functional dependence on several functions of one or two complex variables, and complex integrals. Thus here we present the most general solution for the situation that these functions are identically zero:
\begin{eqnarray*}
F_1&=&\tau_0\tau_1\tau_2\left( C_7\eta_1+C_8\eta_2 \right)+C_{14}\,,\\
F_2&=&-\frac {\left(C_7 \eta_2-C_8\eta_1 \right) \sqrt 
{C_{11}}\tau_2\tau_1\left( 2C_2-\tau_0 \right) \sqrt {C_{11}+C_{12}+C_{13}}}{C_{11}+C_{12}}\\
&&+\frac {\sqrt {C_{13}} \left(2C_4-\tau_1 \right) 
\sqrt {C_{12}} \left( 2C_6-\tau_2 \right) 
 \left( C_7\eta_1+C_8 \eta_2 \right)\tau_0}{C_{11}+C_{12}}+C_{10}\,,\\
F_3&=&\frac {-
\tau_0\left( 2C_4-\tau_1 \right)  \left(C_7 \eta_2-C_8 \eta_1 \right) \sqrt {C_{12}}\sqrt {C_{11}+C_{12}+C_{13}}\tau_2 }{C_{11}+C_{12}}\\
&&-\frac {\sqrt {C_{13}} \left(  C_7\eta_1+C_8\eta_2  \right) \tau_1 \left(2C_6- \tau_2
 \right)  \left(2 C_2-\tau_0 \right) \sqrt {C_{11}} }{C_{11}+C_{12}}+C_9\,,\\
\tau_0&=&\exp(2\sqrt{C_{11}}h_0)C_1+C_2\,,\ 
\tau_1=\exp(2\sqrt{C_{12}}h_1)C_3+C_4\,,  \
\tau_2=\exp(2\sqrt{C_{13}}h_2)C_5+C_6\,,\\
\eta_1&=&\sin( \sqrt {C_{11}+C_{12}+C_{13}}h_3)\exp(-\sqrt{C_{12}}h_1)\exp(-\sqrt{C_{13}}h_2)\exp(-\sqrt{C_{11}}h_0)\,,\\
\eta_2&=&\cos( \sqrt {C_{11}+C_{12}+C_{13}}h_3)\exp(-\sqrt{C_{12}}h_1)\exp(-\sqrt{C_{13}}h_2)\exp(-\sqrt{C_{11}}h_0)\,,
\end{eqnarray*}
where $C_1,\dots,C_{14}$ are constants.
See  {\tt{http://www.math.uni-hamburg.de/home/cortes/CCG.zip}} for the maple worksheet.

\bt\label{theoremFlat}
Consider the quaternionic-skew Hermitian manifold $(M^{4n}, Q, \om)$ $(n>1)$ endowed with the torsion-free quaternionic skew-Hermitian connection $\nabla=\nabla^{Q, \om}$ such that $R^{Q, \om}=0$. 
Consider also the scalar 2-form $ \hat\om+\beta$, with 
\[
\beta=F_1(\dd h_0\wedge \dd h_1+\dd h_2\wedge \dd h_3)+F_2(\dd h_0\wedge \dd h_2-\dd h_1\wedge \dd h_3)+F_3(\dd h_0\wedge \dd h_3+\dd h_1\wedge \dd h_2)\,,
\]
where the  functions $F_a$ are given in Lemma \ref{lemmaFlat}.
Then $\hat\om+\beta$ is closed and the corresponding $\SO^*(2(n+1))$-structure $(\hat{H}, \hat\om+\beta)$ on the Swann bundle $\hat{M}$ is of type $\mc{X}_{57}$. It is a torsion-free $\SO^*(2(n+1))$-structure if and only if the functions $F_a$ are constant for all $a=1, 2, 3$.
\et
\pr
By direct computations, one can check that $\dd \beta=0$ for the above functions $F_a$, and the first claim follows. Clearly, in the exponential coordinates the Obata connection $\nabla^{\hat H}$ is the flat connection on $[\E\Hh]\times \Hn^{\times}/\Z_2$. Thus the $\SO^*(2(n+1))$-structure $(\hat{H}, \hat\om+\beta)$ is torsion-free if and only if $\nabla^{\hat H}(\hat\om+\beta)=0$. Since on $TM$ it coincides with the flat unimodular Oproiu connection and $\hat\om$ is constant along the fibres, we get $\nabla^{\hat H}\hat\om=0$. Clearly,  $\nabla^{\hat H}\beta=\sum_{cycl}^3(\nabla^{\hat H}F_a)\wedge(\dd h_0\wedge \dd h_a+\dd h_b\wedge \dd h_c)$, 
 and thus the $\SO^*(2(n+1))$-structure is torsion-free if and only if $F_a$ are constant for all $a$.
\pro

{ }

\begin{thebibliography}{90}

 \bibitem{AC}
 Alekseevsky, D. V.;  V. Cort\'es.
 \newblock  ``Classification of pseudo-Riemannian symmetric spaces of quaternionic K\"ahler type.''
 \newblock in ``Lie groups and invariant theory.''  Providence, RI: American Mathematical Society 213 (AMS). Translations. Series 2. \textit{Adv. Math. Sci.}  56, 33--62,  (2005).
 
 
  \bibitem{ACDM}
 Alekseevsky, D. V.;   Cort\'es, V.; Dyckmanns, M; T. Mohaupt.
 \newblock ``Quaternionic K\"ahler metrics associated with special K\"ahler manifolds.''
 \newblock \textit{J. Geom. Phys.}, 92, 271--287, (2015).
 
  
\bibitem{AM}
 Alekseevsky, D. V.; S. Marchiafava.
\newblock ``Quaternionic structures on a manifold and subordinated structures.'' 
\newblock \textit{Ann. Mat. Pura Appl.} (IV), Vol  CLXXI, 205-273, (1996).
 
 
 \bibitem{CGWPartI}
  Chrysikos, I.;  Gregorovi\v{c}, J.; H. Winther.
 \newblock  ``Differential geometry of $\SO^\ast(2n)$-type structures.''
 \newblock \textit{Annali di Matematica Pura ed Applicata (1923 -)},   60pp, (2022), (doi.org/10.1007/s10231-022-01212-y). 
  

 \bibitem{CGWPartII}
  Chrysikos, I.;  Gregorovi\v{c}, J.; H. Winther.
  \newblock   ``Differential geometry of $\SO^\ast(2n)$-type structures - Integrability.''
 \newblock  \textit{Analysis and Mathematical Physics}, Vol 12 (93),  1--52,  (2022). 

  
 \bibitem{CH}
Cort\'es, V.; K. Hasegawa.
\newblock   ``The quaternionic/hypercomplex-correspondence.''
 \newblock   \textit{Osaka J. Math.}, 58, 213--238,  (2021).
 

\bibitem{Gil}
Gilmore, R.
 \newblock  ``Lie Groups, Lie Algebras, and Some of Their Applications.''
 \newblock   \textit{A Willey-Interscience Publication},  New-York, 1974.
 
 
 \bibitem{Jan}
 Gregorovi\v{c}, J.
 \newblock ``Geometric structures invariant to symmetries.'' (Phd thesis), 
 \newblock arXiv:1207.0193.
 

 \bibitem{MS1}
 Merkulov, S. A.; L. J. Schwachh\"ofer.
 ``Classification of irreducible holonomies of torsion-free affine connections.''
 \textit{Ann. Math.}, 150, 77--49, (1999). Addendum:
``Classification of irreducible holonomies of torsion-free affine connections.''
 \textit{Ann. Math.}, 150, 1177--1179, (1999).
 

 \bibitem{PPS}
Pedersen, H.;  Poon, Y.;  A. Swann.
 \newblock   ``Hypercomplex structures associated to quaternionic manifolds.''
 \newblock    \textit{Differential Geom. Appl.}, 9, 273--292, (1998).
      
 
\bibitem{Pont}
Pontecorvo, M.
\newblock ``Complex structures on quaternionic manifolds.''
\newblock  \textit{Differ. Geom. Appl.}, 4, 163--177, (1994).
    
   
\bibitem{Salamon86}
Salamon, S.   M.
 \newblock ``Differential geometry of quaternionic manifolds.''
 \newblock \textit{Ann. Scient. Ec. Norm. Sup.}, $4^{e}$ s\'erie,
19,  31--55, (1986).

    
\bibitem{S1} 
 Schwachh\"ofer, J.
\newblock  ``Connections with irreducible holonomy representations.''
 \newblock  \textit{Advances in Mathematics} 160 (1), 1--80, (2001).
 
    
\bibitem{Swann}
Swann, A.
\newblock  ``Hyper-K\"ahler and quaternionic-K\"ahler geometry.''
\newblock  \textit{Math. Ann.}, 289, 421--450,  (1991).


\end{thebibliography}
\end{document}